\documentclass[11pt,a4paper,halfparskip]{article}

\usepackage[hypertex]{hyperref}
\usepackage{bbm}
\usepackage[intlimits,sumlimits]{amsmath}
\usepackage{amsthm}
\usepackage{amscd}
\usepackage{amssymb}
\usepackage[all]{xy}

\newcommand{\cu}{\mathcal{U}}
\newcommand{\co}{\mathcal{O}}
\newcommand{\cm}{\mathcal{M}}

\newcommand{\cv}{\mathcal{V}}
\newcommand{\cs}{\mathcal{S}}
\newcommand{\cf}{\mathcal{F}}
\newcommand{\cl}{\mathcal{L}}
\newcommand{\ca}{\mathcal{A}}

\newcommand{\cp}{\mathcal{P}}

\newcommand{\cg}{\mathcal{G}}

\newcommand{\cj}{\mathcal{J}}

\newcommand{\fr}{\mathfrak{R}}

\newcommand{\cgl}{\mathfrak{GL}}

\newcommand{\End}{\mathrm{End}}

\newcommand{\id}{\mathrm{id}}

\newcommand{\cem}{\mathrm{cem}}

\newcommand{\Tr}{\mathrm{Tr}}

\newcommand{\sdiff}{\widehat{\mathcal{SD}}}

\newcommand{\realR}{\mathbbm{R}}
\newcommand{\complexC}{\mathbbm{C}}
\newcommand{\intZ}{\mathbbm{Z}}
\newcommand{\fieldK}{\mathbbm{K}}
\newcommand{\naturalN}{\mathbbm{N}}

\newcommand{\olv}{\overline{V}}
\newcommand{\olw}{\overline{W}}
\newcommand{\olr}{\overline{\realR}}
\newcommand{\ovlr}[1]{\overline{\realR^{#1}}}
\newcommand{\olc}{\overline{\complexC}}
\newcommand{\olk}{\overline{\fieldK}}

\newcommand{\ovlk}[1]{\overline{\fieldK^{#1}}}
\newcommand{\olf}{\overline{f}}
\newcommand{\olcdot}{\overline{\cdot}}
\newcommand{\olpi}{\overline{\Pi}}
\newcommand{\olmu}{\overline{\mu}}

\newcommand{\Hom}{\mathrm{Hom}}
\newcommand{\ihom}{\underline{\mathrm{Hom}}}

\newcommand{\iend}{\underline{\mathrm{End}}}

\newcommand{\Spec}{\mathrm{Spec}\,}
\newcommand{\Aut}{\mathrm{Aut}}

\newcommand{\Sym}{\mathrm{Sym}}

\newcommand{\catgr}{\mathsf{Gr}}
\newcommand{\catsets}{\mathsf{Sets}}
\newcommand{\catc}{\mathsf{C}}
\newcommand{\catd}{\mathsf{D}}

\newcommand{\catmod}{\mathsf{Mod}}
\newcommand{\catsmod}{\mathsf{SMod}}
\newcommand{\catsrepmod}{\mathsf{SRepMod}}
\newcommand{\cattop}{\mathsf{Top}}
\newcommand{\catman}{\mathsf{Man}}
\newcommand{\catbsdom}{\mathsf{BSDom}}
\newcommand{\catbsman}{\mathsf{BSMan}}
\newcommand{\catsman}{\mathsf{SMan}}

\newcommand{\catsvect}{\mathsf{SVec}}
\newcommand{\catvect}{\mathsf{Vec}}
\newcommand{\catspoint}{\mathsf{SPoint}}

\newcommand{\catfinsman}{\mathsf{FinSMan}}
\newcommand{\catsalg}{\mathsf{SAlg}}
\newcommand{\catalg}{\mathsf{Alg}}
\newcommand{\catslie}{\mathsf{SLie}}
\newcommand{\catlie}{\mathsf{Lie}}

\newtheorem{thm}{Theorem}[section]
\newtheorem{prop}[thm]{Proposition}
\newtheorem{lemma}[thm]{Lemma}
\newtheorem{cor}[thm]{Corollary}
\newtheorem{dfn}[thm]{Definition}

\newtheorem{rmk}[thm]{Remark}

\begin{document}

\author{Christoph Sachse\\[0.5em]
Max Planck Institute for Mathematics in the Sciences,\\ Inselstr. 22, D-04103 Leipzig, Germany}
\title{A Categorical Formulation of Superalgebra and Supergeometry}
\date{}

\maketitle

\begin{abstract}
We reformulate superalgebra and supergeometry in completely categorical terms by a
consequent use of the functor of points. The increased abstraction of this approach 
is rewarded by a number of great advantages. First, we show that one can extend supergeometry completely
naturally to infinite-dimensional contexts. Secondly, some subtle and sometimes obscure-seeming points
of supergeometry become clear in light of these results, e.g., the relation between the
Berezin-Leites-Kostant and the de Witt-Rogers approaches, and the precise geometric meaning of odd 
parameters in supergeometry. In addition, this method allows us to construct in an easy manner superspaces
of morphisms between superobjects, i.e., the inner Hom objects associated with the sets of such morphisms.
The results of our work rely heavily and were inspired by the ideas of
V. Molotkov, who in \cite{M:Infinite} outlined the approach extensively expounded here.
\end{abstract}

\section{Introduction}

Superalgebra and supergeometry are fields of mathematics which owe their invention to a physical
concept: supersymmetry. In quantum field theory, fields of integer spin (called bosons) are described
by commuting variables while those of half-integer spin (fermions) have to be anticommuting in
order for the theory to produce meaningful results. It was Berezin who had the idea to combine
these physically very differently behaved types of fields into a single $\intZ_2$-graded space.
Operators, as well, could then form $\intZ_2$-graded algebras, and their operation would have to
preserve parity. The radically new option that this formalism introduces is that one can also
write down transformations which are proportional to odd parameters. This made it possible to
transform bosons into fermions and vice versa by means of fermionic parameters. Supersymmetry is then
nothing else than the statement that a physical theory remains invariant under this type of
transformations.

In mathematical terms,  this means that supersymmetry does not simply consist of the idea to
use $\intZ_2$-graded
spaces and algebras, but also to carry out all operations over a $\intZ_2$-graded base.
Physical applications, in particular in gravity theory, very soon made it desirable to extend this concept from pure algebra to
geometry. This geometry would have to include odd as well as even coordinates and treat them
on a completely equal footing. It was found out already in the seventies that this is indeed 
possible \cite{L:Introduction}.
Again, the full meaning of supergeometry only becomes visible
by working over a base which is itself ``super'' in order to allow for odd parameters in all
coordinate changes. The resulting theory extends commutative to ``supercommutative'' geometry.
Its noncommutativity is so tame (an actual noncommutative geometer would probably not even call it
noncommutative) that it still allows to carry over almost all objects of classical
geometry, but it also already exhibits completely new phenomena.

The correct description of supermanifolds required several technical tools
(e.g., sheaf theory, point functors, working in families, etc.) which were new to physics at that
time, and this fact has made the relationship between physicists and mathematicians working on
super questions at times a difficult one. From a mathematical point of view, the ringed space approach
common in algebraic geometry is the obvious and essentially the only way to extend the ideas
of superalgebra to the realm of geometry. Rings are the starting point for the construction of any
geometric space in algebraic geometry anyway, and
the transition from commutative to supercommutative ones can be handled in a completely natural and very
elegant way in this formalism. The price, however, is a considerable abstraction and 
technical complexity.

Therefore, various other approaches have been advocated, most prominently the Rogers-de Witt approach.
These approaches mostly tried to avoid the technicalities of sheaves and point functors and tried
to reproduce a closer analogue of ordinary differential geometry. Many practicioners working 
on supersymmetric field theories and
string theory have found these approaches easier to handle, and they were sometimes claimed to be
``more adequate to physics'' than the algebraic approach. The categorical formulation outlined below will
allow us make the relation between these approaches completely transparent.
In particular, we will argue that the more concrete looking de Witt approach
is indeed correct as long as one constructs everything functorially in the Grassmann algebra that one
chooses (and infinite-dimensional Grassmann algebras are unnecessary). Although this fact seems to have been
rarely spelled out in the literature, this is the way the
Rogers-de Witt geometry has been used in most applications.
In this case, it is equivalent to the use of just one set of points of the corresponding point
functor of the super object under study. Any kind of topologization of the odd dimensions of a supermanifold,
however, will produce faulty results, since it destroys functoriality with respect to base change.

The categorical approach should
not be understood as simply one more way to phrase the same old results, but rather as a complement 
to the ringed space approach.
It is based on the exploitation of a tool which is standard  in algebraic geometry, namely
the functor of points supplied by the Yoneda Lemma, which describes an object of a given
category by the sets of morphisms of other objects of that category into it. The usefulness of this
method was, of course, well known to the founders of the ringed space approach \cite{L:Introduction},
\cite{BL:Supermanifolds}, \cite{M:Gauge}, and was used, e.g., for the definition of supergroups 
\cite{SoS}. From its definition alone it is, however, often completely unclear how one could get
a manageable description of this functor which allows concrete calculations and constructions.
V. Molotkov in \cite{M:Infinite} outlined a program that reformulated superalgebra as
well as supergeometry entirely in this categorical framework and also proposed a way to make
the functor of points really usable, namely by presenting a countable set of generators for
the category of supermanifolds. Unfortunately, his preprint \cite{M:Infinite}
contained no proofs, and is not easily accessible for more practically or physically oriented researchers. 
In this work we carry out most of the program proposed by Molotkov and develop supermathematics from
a categorical perspective. Similar ideas have been investigated in \cite{S:definition}, \cite{V:Maps}.

Our original motivation to investigate these questions was the need for
a practicable model for infinite-dimensional supergeometry. Ringed spaces become inadequate in
infinte dimensions, and one has to resort to concepts like functored spaces. This is not a super
problem, but occurs as well in ordinary infinite-dimensional geometry \cite{D:probleme}. The categorical
construction solves this problem for supermanifolds. As an example it will be used in a subsequent paper 
\cite{S:Structure} to explicitly construct the diffeomorphism supergroup of a supermanifold
$\cm$.

\subsection*{Acknowledgments}

I thank the Klaus Tschira Stiftung and the MPI for Mathematics in the Sciences for financial support.
The results of this work were in large parts inspired and some already announced in Vladimir Molotkov's 
preprint \cite{M:Infinite}. Besides these ideas I owe him thanks for continued advice and corrections
while trying to forge this bunch of ideas into a neat piece of mathematics. I am also grateful to
J\"urgen Jost, Dimitri Leites, Guy Buss and Brian Clarke for encouragement, advice and valuable criticism.

\section{Brief review of superalgebra and supergeometry}

This section contains a very brief review of the key ideas and constructions of superalgebra and supergeometry. 
For in-depth (mathematical) treatments of the matters 
in this section see, e.g., \cite{DM:Notes}, \cite{SoS}, \cite{V:Supersymmetry}.

\subsection{Linear superspaces and superalgebras}

The elements of $\intZ_2=\intZ/2\intZ$ will be denoted as $\{\bar{0},\bar{1}\}$. The field $\fieldK$
denotes either $\realR$ or $\complexC$. We begin by recalling that a \emph{superring} $R$ is simply
a $\intZ_2$-graded ring, i.e. $R=R_{\bar{0}}\oplus R_{\bar{1}}$, and that a morphism of superrings is
a morphism of graded rings.

\begin{dfn}
\label{def:sring}
A module $M$ over a superring $R$ is called a supermodule, if it is $\intZ_2$-graded,
\begin{equation}
M=M_{\bar{0}}\oplus M_{\bar{1}},
\end{equation}
and if $R_{\bar{i}}\cdot M_{\bar{j}}\subseteq M_{\bar{i}+\bar{j}}$.
The submodule $M_{\bar{0}}$ is called even, $M_{\bar{1}}$ 
is called odd. A morphism $\phi:M\to M'$ of
$R$-supermodules is a morphism of $R$-modules which preserves the grading.
\end{dfn}

We denote by $\catsmod_R$ (resp. ${}_R\catsmod$) the category of right (resp. left) $R$-supermodules.
Note that any ring $R$ can be considered as a superring by just setting $R_{\bar{0}}:=R$ and
$R_{\bar{1}}:={0}$. 

An element $m\in M_{\bar{i}}$ is said to be homogeneous of parity $p(m)=i$.
An inhomogeneous element is said to be of indefinite parity.

A supermodule over a field $\fieldK$ is called a $\fieldK$-super vector space. 
The category of $\fieldK$-super vector spaces will be denoted by $\catsvect_\fieldK$. Particularly important
are the standard super vector spaces $\fieldK^{m|n}$, which have $m$ even and $n$ odd dimensions.

\begin{dfn}
A $\fieldK$-superalgebra $A$ is a $\fieldK$-super vector space endowed with a $\fieldK$-bilinear morphism
\begin{equation}
\mu:A\times A\to A.
\end{equation}
The algebra
$A$ is called supercommutative if for all homogeneous elements $a,b$, one has
\begin{equation}
\label{salgcomm}
\mu(a,b)=(-1)^{p(a)p(b)}\mu(b,a).
\end{equation}
$A$ is associative, resp. with unit, if it is associative, resp. has a unit, as an ordinary $\fieldK$-algebra.
\end{dfn}

Lie superalgebras are defined completely analogously.
Expressions like (\ref{salgcomm}) are extended to inhomogeneous elements by linearity.
It is clear that if a unit exists, it has to be even. If no confusion can arise, 
multiplication in an algebra will be
denoted by either a dot (like $a\cdot b$) or by juxtaposition (like $ab$). In the following, 
we will assume that all superalgebras are associative and have a unit.

For a supercommutative algebra, every left module can be made a right module by defining
\begin{equation}
\label{modstr}
m\cdot a:=(-1)^{p(a)p(m)}a\cdot m
\end{equation}
for all homogeneous elements $a\in A$ and $m\in M$, and extending this definition by linearity. From now on,
we will restrict ourselves to supercommutative algebras with unit. Consequently, we will just speak of
``modules'', all of which are defined to be left-modules whose right module structure is given by
(\ref{modstr}).

Direct sums, tensor products and $\Hom$-spaces can be defined in an obvious way for supermodules. Let
$A$ be a supercommutative $\fieldK$-superalgebra and $M,N$ be $A$-modules. Then
\begin{eqnarray}
(M\oplus N)_{\bar{i}} &:=& M_{\bar{i}}\oplus N_{\bar{i}},\\
(M\otimes_A N) &:=& (M\otimes_{\fieldK} N)/I,
\end{eqnarray}
where $I$ is the ideal generated by elements of the form $ma\otimes n-m\otimes an$, and
\begin{equation}
(M\otimes_{\fieldK} N)_{\bar{i}}:=\bigoplus_{\bar{j}+\bar{k}=\bar{i}}M_{\bar{j}}\otimes_\fieldK N_{\bar{k}}.
\end{equation}
Finally one habitually sets
\begin{equation}
\ihom_A(M,N)_{\bar{i}} := \{f:M\to N\mid f(M_{\bar{j}})\subset N_{\bar{j}+\bar{i}}\}
\end{equation}
where on the right hand side, $f$ is meant to be an ordinary module homomorphism between $M$ and $N$.
It is important to note here that $\Hom_A(M,N)$ does in general \emph{not} coincide with $\ihom_A(M,N)$: the
former is the set of morphisms of supermodules, which by definition preserve parity. This set does not carry
any natural supermodule structure itself. The latter is, strictly speaking, the inner Hom object associated
to this set. One always has $\ihom_A(M,N)_{\bar{0}}=\Hom_A(M,N)$, that is, the inner Hom object extends the
set of morphisms of supermodules by adding parity reversing morphisms. The distinction between the two objects
may seem overly picky at this point, but it will turn out that inner Hom objects play a particularly important
role in supermathematics. In supergeometry, e.g., as diffeomorphism supergroups, their construction is by 
far less obvious than in superalgebra.

The above definition of the tensor product together with the following choice of the commutativity isomorphisms
\begin{eqnarray}
\label{tprodcomiso}
c_{V,W}:V\otimes_A W &\to& W\otimes_A V\\
\nonumber
v\otimes w &\mapsto& (-1)^{p(v)p(w)}w\otimes v.
\end{eqnarray}
make $\catsmod_A$ a braided tensor category. This choice of the braiding encodes the so-called sign rule
which asserts that whenever two neighbouring factors in a multiplicative expression involving elements of supermodules are
interchanged, one picks up a factor $(-1)$ and thus is the key difference between superalgebra and ordinary
algebra.

The braiding isomorphisms supply the notions of symmetric and exterior powers as well as those of
symmetric and antisymmetric maps in the usual fashion: one sets $\Sym^n(V)$ to be the quotient of
$V^{\otimes n}$ divided by the action of the permutation group $S_n$ described by $c_{V,V}$. Analogously,
$\wedge^n(V)$ is defined as $V^{\otimes n}$ divided by the action of $S_n$ times its sign character. Clearly,
for $V=V_{\bar{0}}\oplus V_{\bar{1}}$,
\begin{eqnarray}
\Sym^n(V) &=& \bigoplus_{k=0}^n\Sym^k(V_{\bar{0}})\otimes\wedge^{n-k}(V_{\bar{1}}),
\end{eqnarray}
where $\Sym$ and $\wedge$ on the right hand side denote the ordinary operations on vector spaces.
Symmetric (resp. ``supersymmetric'') multilinear maps are then
morphisms $f:V^{\otimes n}\to W$ which are invariant under the action of $S_n$ by $c_{V,V}$. In other words,
they descend to well defined maps $f:\Sym^n(V)\to W$. We will denote the set of supersymmetric morphisms
between $V^{\otimes n}$ and $W$ by $\Sym^n(V;W)$. The sets $\Sym^n(V;W)$ are again the even subspaces
of inner Hom objects $\underline{\Sym}^n(V;W)$.

\begin{dfn}
\label{def:copf}
Let $R$ be a superring. The change of parity functor is the functor $\Pi:\catsmod_R\to\catsmod_R$ 
(equivalently for left modules) which assigns to a supermodule 
$M=M_{\bar{0}}\oplus M_{\bar{1}}$ the supermodule $\Pi(M)$ with $(\Pi(M))_{\bar{0}}=M_{\bar{1}}$
and $(\Pi(M))_{\bar{1}}=M_{\bar{0}}$.
\end{dfn}

Parity reversal has to be a functor since any morphism between two supermodules has to preserve
parity.
One clearly has $\Pi(\fieldK)=\fieldK^{0|1}$ and for any $\fieldK$-super vector space $V$,
\begin{equation}
\fieldK^{0|1}\otimes_\fieldK V\cong\Pi(V).
\end{equation}

\subsubsection{The category of finite-dimensional Grassmann algebras}
\label{subsect:gr}

A particularly important role will be played by the finitely generated Grassmann algebras $\Lambda_n$, i.e.,
the free commutative superalgebras on $n$ odd generators. In particular, we will use the category $\catgr$, which
contains exactly one object $\Lambda_n$ for each integer $n\geq 0$ (i.e., the category of isomorphism classes of Grassmann algebras) 
and whose morphisms are the superalgebra morphisms of the $\Lambda_n$.

Note that the field $\fieldK=\Lambda_0$ is the null object of $\catgr$: the terminal morphisms $\epsilon_{\Lambda_n}$
act by removing all odd generators, while the initial morphisms $c_{\Lambda_n}$ inject $\fieldK$ into $\Lambda_n$.

We will write $\catgr$ and $\Lambda_n$ both for the real and complex Grassmann algebras 
and only distinguish them when the ground field really matters. 

\begin{prop}
\label{homgr}
Let $\Lambda_n,\Lambda_m$ be the Grassmann algebras over $\fieldK$ on $n$ and $m$ generators, respectively.
Then there exists an isomorphism of $\fieldK$-vector spaces
\begin{equation}
\Hom_\catgr(\Lambda_n,\Lambda_m)\cong \fieldK^n\otimes_\fieldK \Lambda_{m,\bar{1}}.
\end{equation}
\end{prop}
\begin{proof}
Let $\xi_1,\ldots,\xi_n$ be the free generators of $\Lambda_n$.
A morphism $\phi:\Lambda_n\to\Lambda_m$ is a homomorphism of $\fieldK$-algebras which preserves parity.
Thus $\phi(1)=1$, and $\phi$ is uniquely determined by choosing the images $\phi(\xi_n)$ of its
generators which have to lie in $\Lambda_{m,\bar{1}}$. Setting 
$V=\mathrm{Span}_\fieldK(\xi_1,\ldots,\xi_n)$, we can write
\[
\Hom_\catgr(\Lambda_n,\Lambda_m)\cong \Lambda_{m,\bar{1}}\otimes_\fieldK V^*,
\]
where $V^*$ is the dual space of $V$.
\end{proof}

\subsection{Supermanifolds}
\label{sect:sman}

Here, again, we will only very briefly review the key ideas of supergeometry. See \cite{SoS}, \cite{V:Supersymmetry},
\cite{DM:Notes} for more complete treatments. 

Supergeometry is based on the idea to replace the commutative rings of ``functions'' which make up the structure
sheaves of commutative geometric spaces by supercommutative ones. In contrast to general noncommutative geometry,
this works out within the classical framework of algebraic geometry, because almost all crucial tools,
like, e.g., localization, are amenable to supercommutative rings.

\begin{dfn}
A locally superringed space $\cm=(M,\co_\cm)$ is a topological space $M$ endowed with a sheaf $\co_\cm$ 
of local supercommutative rings.
A morphism of superspaces is a morphism of locally ringed spaces which is stalkwise a homomorphism of
supercommutative rings.
\end{dfn}

The structure sheaf $\co_\cm$ contains a canonical subsheaf of nilpotent ideals 
generated by odd elements:
\begin{equation}
\label{nilideal}
\cj = \co_{\cm,\bar{1}}\oplus\co_{\cm,\bar{1}}^2.
\end{equation}
The superspace $\cm_{red}=(M,\co_\cm/\cj)$ is then purely even and possesses a canonical embedding 
\begin{equation}
\cem:\cm_{red}\to\cm.
\end{equation}
The space $\cm_{red}$ is called the underlying space of $\cm$.
In general, $\cm_{red}$ can itself be a non-reduced space. The ``completely reduced'' space is
denoted as $\cm_{rd}$. In this work we will not consider non-reduced underlying spaces,
so we always assume $\cm_{red}=\cm_{rd}$.

A particular role is played by the model spaces $\realR^{m|n}$ and $\complexC^{m|n}$, which are
defined by
\begin{eqnarray}
\realR^{m|n} &:=& (\realR^m,C^\infty_{\realR^m}[\theta_1,\ldots,\theta_n])\\
\complexC^{m|n} &:=& (\complexC^m,\co_{\complexC^m}[\theta_1,\ldots,\theta_n])
\end{eqnarray}
where $C^\infty$ and $\co$ denote the sheaves of smooth and holomorphic functions, respectively, and
$\theta_1,\ldots,\theta_n$ are odd generators, which freely generate the structure sheaf of
these superspaces over $C^\infty$ and $\co$, respectively.

\begin{dfn}
A real, resp. complex, supermanifold of dimension $m|n$ is a superringed space which is locally isomorphic
to $\realR^{m|n}$, resp. $\complexC^{m|n}$.
\end{dfn}

Thus, strictly speaking, a super vector space is \emph{not} a supermanifold. Each super vector space $V$ gives
rise to a supermanifold $\olv=(V_{\bar{0}},C^\infty_{V_{\bar{0}}}\otimes\wedge(V_{\bar{1}}))$ (here
$\wedge(V_{\bar{1}})$ denotes the ordinary exterior algebra over $V_{\bar{1}}$), but these two objects behave
rather differently. For example, the underlying set-theoretical model of $V$ is a $(m+n)$-dimensional space, while that of
$\olv$ is only a $m$-dimensional one.

Together with their morphisms as superringed spaces, smooth finite-di\-men\-sion\-al supermanifolds form a 
category which we will
denote as $\catfinsman$.

Sections of the structure sheaf $\co_\cm$ of a supermanifold $\cm$ are called superfunctions. Each superfunction
can locally be expanded into powers of the odd generators:
\begin{equation}
F(x_1,\ldots,x_m,\theta_1,\ldots,\theta_n)=\sum_{I\subseteq\{1,\ldots,n\}}f_I(x_1,\ldots,x_m)\theta_I,
\end{equation}
where the sum runs over all increasingly ordered subsets, $\theta_I$ is the product of the appropriate local odd coordinates,
and the $f_I$ are ordinary smooth (resp. holomorphic) functions. The value of $F$ at $x$ is defined to be $f_0(x)$,
and obviously does not determine $F$, in contrast to ordinary geometry.

\subsubsection{The category of superpoints}
\label{sect:spoint}

\begin{dfn}
A finite-di\-men\-sion\-al supermanifold whose underlying manifold is a one-point topological
space is called a superpoint.
\end{dfn}

Superpoints will play an important role later on, analogous to that of the spaces $\Spec\fieldK$ for
ordinary geometry.
From the above discussion it is evident that superpoints are the supermanifolds associated to
purely odd super vector spaces, i.e., to those of dimension $0|n$. 
Together with their morphisms as supermanifolds, finite-dimensional superpoints form a category $\catspoint$
which is a full subcategory of $\catsman$.

\begin{prop}[see also \cite{M:Infinite}]
\label{spointeq}
There exists an equivalence of categories
\begin{equation}
\cp:\catgr^\circ\to\catspoint.
\end{equation}
\end{prop}
\begin{proof}
Define the functor $\cp$ on the objects $\Lambda_n$ of $\catgr^\circ$ as
\begin{equation}
\cp(\Lambda_n):=\cp_n:=\Spec\Lambda_n=(\{*\},\Lambda_n).
\end{equation}
To every morphism $\varphi:\Lambda_n\to\Lambda_m$ of Grassmann algebras, assign the morphism
\begin{equation}
\Phi=(\id_{\{*\}},\varphi):\cp_m \to \cp_n.
\end{equation}
To see that this functor establishes an equivalence, note first that it is fully faithful: on the set
of morphisms it is a bijection from $\Hom_\catgr(\Lambda_n,\Lambda_m)$ to $\Hom_\catsman(\cp_m,\cp_n)$.
The last property to check is essential surjectivity, i.e., every superpoint has to be isomorphic to one
of the $\cp_n$. This is clear from the fact that since a superpoint is a supermanifold its structure
sheaf must be a free superalgebra on $n$ odd generators. Since any such algebra is isomorphic to
$\Lambda_n$, the assertion is proved.
\end{proof}

This allows us to restrict our attention to the supermanifolds $\Spec\Lambda_n$ when talking about superpoints.

\section{Superalgebra in the categorical setting}

As a first step we will translate linear and commutative
superalgebra into the language of the functor of points. This construction
can be viewed as a systematic treatment of the so-called even rules principle, which
is a way to do superalgebra without having to handle odd quantities \cite{DM:Notes}, \cite{V:Supersymmetry}.

We will use the following notational conventions throughout the rest of this paper: if $C$ is an object
of a category $\catc$, we write simply $C\in\catc$. The category of functors $\catc\to\catd$ between
the small category $\catc$ and an arbitrary category $\catd$ will be denoted as $\catd^\catc$. 
For in-depth accounts of the categorical tools used below, see \cite{J:Sketches}, 
\cite{M:Categories}, \cite{MM:Sheaves}.

\subsection{The rings $\olr$ and $\olc$}
In the following, we will almost exclusively be concerned with the functor categories
\[
\catsets^\catgr\quad\supset\quad\cattop^\catgr\quad\supset\quad\catman^\catgr,
\]
that is, covariant functors from the category $\catgr$ of Grassmann algebras defined in \ref{subsect:gr} into
the category of sets, topological spaces and smooth Banach manifolds, respectively.

For $\fieldK=\realR,\complexC$, define a functor $\olk\in\catsets^{\catgr}$ by
\begin{eqnarray}
\olk(\Lambda) &:=& (\Lambda\otimes_{\fieldK} \fieldK)_{\bar{0}}=\Lambda_{\bar{0}},\\
\olk(\varphi) &:=& \varphi\big|_{\Lambda_{\bar{0}}}
\end{eqnarray}
for $\varphi:\Lambda\to{\Lambda'}$ a morphism in $\catgr$,
the category of Grassmann algebras over the field $\fieldK$. Clearly, $\olk$ is a commutative ring with unit
in $\catsets^{\catgr}$: each $\olk(\Lambda)$ has a commutative ring structure inherited from $\Lambda_{\bar{0}}$,
with a unit induced by the unique morphism $\Lambda\to\fieldK$.

The rings $\olk$ will replace the ground field $\fieldK$ in the categorical definition of super spaces.

\subsection{$\olk$-modules in $\catsets^\catgr$}

In this section we will introduce a particular class of $\olk$-modules in
$\catsets^\catgr$, which are the avatars super vector spaces over $\fieldK$ in
the categorical formulation.

Let $V$ be some real or complex super vector space. We define a functor $\olv\in\catsets^\catgr$ by
setting
\begin{eqnarray}
\label{def:olv}
\olv(\Lambda) &:=& (\Lambda\otimes_\fieldK V)_{\bar{0}}=\left(\Lambda_{\bar{0}}\otimes V_{\bar{0}}\right)\oplus%
\left(\Lambda_{\bar{1}}\otimes V_{\bar{1}}\right)\\
\nonumber
\olv(\varphi) &:=& (\varphi\otimes\mathrm{id}_{V})\big|_{\olv(\Lambda)}\qquad\mathrm{for}\quad%
\varphi:\Lambda\to\Lambda'.
\end{eqnarray}
All sets $\olv(\Lambda)$ are naturally $\Lambda_{\bar{0}}$-modules, and thus, $\olv$ is a $\olk$-module.

Let $f:V_1\times\ldots\times V_n\to V$ be a multilinear map of $\fieldK$-super vector spaces. To $f$,
we assign the functor morphism $\olf:\olv_1\times\ldots\times\olv_n\to\olv_n$ whose components
\begin{equation}
\olf_\Lambda:\olv_1(\Lambda)\times\ldots\times\olv_n(\Lambda)\to\olv(\Lambda)
\end{equation}
are defined by
\begin{equation}
\label{fbar}
\olf_\Lambda(\lambda_1\otimes v_1,\ldots,\lambda_n\otimes v_n):=\lambda_n\cdots\lambda_1\otimes%
f(v_1,\ldots,v_n)
\end{equation}
for all $\lambda_i\otimes v_i\in\olv_i(\Lambda)$. All maps $\olf_\Lambda$ are clearly $\Lambda_{\bar{0}}$-linear,
hence $\olf$ is a $\olk$-$n$-linear morphism in $\catmod_{\olk}(\catsets^\catgr)$.

Given some commutative ring $R$, we will denote the set of $R$-$n$-linear morphisms
$M_1\times\ldots\times M_n\to M$ by $L^n_R(M_1,\ldots,M_n;M)$.
$L^n_R$ is naturally an $R$-module itself, and if the $M_i,M$ are supermodules,
$L^n_R$ it is identical to the inner Hom-object, i.e., a supermodule itself. 

So, in particular, $L^n_{\olk}$ is a $\olk$-module which implies that it is a
$\fieldK$-module, because $\Lambda_{\bar{0}}$-linearity always entails $\fieldK$-linearity. It turns out that
the $\fieldK$-module structures on $L^n_{\fieldK,\bar{0}}$ and $L^n_{\olk}$ coincide as a consequence of functoriality.

\begin{prop}[see also \cite{M:Infinite}\footnote{I am grateful to V. Molotkov for pointing out an 
error in my original proof of this
Proposition and for sending me a correct version.}]
\label{isokmod}
The assignment $f\mapsto\olf$, which is a map
\begin{equation}
L^n_{\fieldK,\bar{0}}(V_1,\ldots,V_n;V)\to L^n_{\olk}(\olv_1,\ldots,\olv_n;\olv)
\end{equation}
for any tuple $V_1,\ldots,V_n,V$ of $\fieldK$-super vector spaces, is an isomorphism of $\fieldK$-modules.
\end{prop}
\begin{proof}
Definition (\ref{fbar}) assigns to every $f\in L^n_{\fieldK,\bar{0}}(V_1,\ldots,V_n;V)$ a functor
morphism $\olf$. We have to
show that one can reconstruct a unique $f$ from a given $\olf\in L^n_{\olk}(\olv_1,\ldots,\olv_n;\olv)$. 
Let a seqeuence $\lambda_i\otimes v_i\in \olv_i(\Lambda)$
be given, $1\leq i\leq n$ and let $j\leq n$ of these $v_i$ be odd, assuming for simplicity that these are
the first $j$. By $\Lambda_{\bar{0}}$-linearity we then have
\begin{equation}
\olf_\Lambda(\lambda_1\otimes v_1,\ldots,\lambda_n\otimes v_n)=%
\lambda_n\cdots\lambda_{j+1}\olf_\Lambda(\lambda_1\otimes v_1,\ldots,\lambda_j\otimes v_j,%
1\otimes v_{j+1},\ldots,1\otimes v_n).
\end{equation}
Now consider the morphism $\eta:\Lambda_j\to\Lambda$ which is defined by $\eta(\theta_k)=\lambda_k$,
where $\theta_k$, $1\leq k\leq j$ are the odd generators of $\Lambda_j$. In order to prove the
statement of the Proposition it will now be enough to show that there exists a unique
\[
g\in L^n_{\fieldK,\bar{0}}(V_1,\ldots,V_n;V)
\]
such that
\begin{equation}
\label{unig}
\olf_{\Lambda_j}(\theta_1\otimes v_1,\ldots,\theta_j\otimes v_j,1\otimes v_{j+1},\ldots,1\otimes v_n)=%
\theta_j\cdots\theta_1\otimes g(v_1,\ldots,v_n).
\end{equation}
This is sufficient because the fact that $f$ is a functor morphism implies that we have a commutative
square
\begin{equation}
\label{lambfunc}
\xymatrix{
\olv_1(\Lambda_j)\times\ldots\times\olv_n(\Lambda_j) \ar[rrr]^{\olv_1(\eta)\times\ldots\times\olv_n(\eta)}
\ar[d]_{\olf_{\Lambda_j}} &&& \olv_1(\Lambda)\times\ldots\times\olv_n(\Lambda) \ar[d]^{\olf_\Lambda} \\
\olv(\Lambda_j) \ar[rrr]^{\olv(\eta)} &&& \olv(\Lambda)
}
\end{equation}
and this then entails that
\[
\olf_\Lambda(\lambda_1\otimes v_1,\ldots,\lambda_n\otimes v_n)=%
\lambda_n\cdots\lambda_1\otimes g(v_1,\ldots,v_n)
\]
for a unique even linear map $g$.

To prove that a unique $g$ as claimed in eq. (\ref{unig}) exists, we first observe that the most general
expression that could appear on the right hand side of eq. (\ref{unig}) reads
\begin{multline}
\label{gensum}
f_{\Lambda_j}(\theta_1\otimes v_1,\ldots,\theta_j\otimes v_j,1\otimes v_{j+1},\ldots,1\otimes v_n)=%
\theta_j\cdots\theta_1\otimes g(v_1,\ldots,v_n)+\\
+\sum_{\substack{m<j \\ j_1>\ldots> j_m}}%
\theta_{j_1}\cdots\theta_{j_m}\otimes g_{j_1\cdots j_m}(v_1,\ldots,v_n),
\end{multline}
where all $g_{j_1\cdots j_m}$ are linear maps. To show that the sum on the right hand side equals zero, 
we again use functoriality (\ref{lambfunc}), this time for the morphisms 
\begin{eqnarray*}
\varphi_l:\Lambda_j &\to& \Lambda_j\\
\varphi_l(\theta_k) &=& 0\qquad\mathrm{if}\quad k=l\\
\varphi_l(\theta_k) &=& \theta_k\qquad\mathrm{if}\quad k\neq l.
\end{eqnarray*}
This evidently kills all summands which contain $\theta_l$ and yields
\[
0=\sum_{\substack{m<j \\ l\notin\{j_1,\ldots, j_m\}}}%
\theta_{j_1}\cdots\theta_{j_m}\otimes g_{j_1\cdots j_m}(v_1,\ldots,v_n).
\]
We obtain such an equation for each $1\leq l\leq j$, therefore the sum on the right hand side of
(\ref{gensum}) must equal zero.
\end{proof}

This Proposition is of central importance for all further constructions. As a first application, we obtain
the following Corollary.

\begin{cor}[see also \cite{M:Infinite}]
\label{cor:barfunct}
The assignment $V\mapsto\olv$ and $f\mapsto\olf$ defines a fully faithful functor
\begin{equation}
\olcdot:\catsmod_\fieldK(\catsets)\to\catmod_{\olk}(\catsets^\catgr).
\end{equation}
\end{cor}
\begin{proof}
It has to be shown that the assignment $f\mapsto\olf$ is a bijection
\begin{equation}
(L_{\fieldK})_{\bar{0}}(V;V')\to L^1_{\olk}(\olv;\overline{V'})
\end{equation}
for any pair $V,V'$ of $\fieldK$-super vector spaces. This is a special case of Prop.~\ref{isokmod}
if one inserts there $V_1=V$ and $V=V'$.
\end{proof}

Actually one gets even more, since Prop.~\ref{isokmod} makes a statement about multilinear 
maps on an arbitrary finite number of arguments.

\begin{cor}[see also \cite{M:Infinite}]
\label{barfunct}
The functor $\olcdot:\catsmod_\fieldK(\catsets)\to\catmod_{\olk}(\catsets^\catgr)$ induces fully
faithful functors
\begin{eqnarray}
\olcdot:\catslie_\fieldK(\catsets) &\to& \catlie_{\olk}(\catsets^\catgr)\\
\olcdot:\catsalg_\fieldK(\catsets) &\to& \catalg_{\olk}(\catsets^\catgr)
\end{eqnarray}
between the categories of $\fieldK$-super Lie algebras and ordinary $\olk$-Lie algebras in $\catsets^\catgr$,
and between $\fieldK$-super algebras and ordinary $\olk$-algebras.
\end{cor}
\begin{proof}
Clearly, $\olcdot$ maps a morphism of $\fieldK$-super algebras into a unique morphism of the corresponding $\olk$-algebras.
We have to check that this is also a surjective assignment.
Let $\olf:(\overline{A},\overline{\mu})\to(\overline{A'},\overline{\mu'})$ be a morphism
of $\olk$-algebras, where $\overline{\mu},\overline{\mu'}$ are the multiplications. Then for every $\Lambda\in\catgr$,
we have
\begin{equation}
\olf_\Lambda(\overline{\mu}_\Lambda(\lambda_1\otimes a_1,\lambda_2\otimes a_2))=%
\overline{\mu}'_\Lambda(\olf_\Lambda(\lambda_1\otimes a_1),\olf(\lambda_2\otimes a_2)).
\end{equation}
for all $\lambda_i\otimes a_i\in \overline{A}(\Lambda)$. The left hand side is
\begin{equation}
\label{lab1}
\olf_\Lambda(\lambda_2\lambda_1\otimes\mu(a_1,a_2))=\lambda_2\lambda_1\otimes f(\mu(a_1,a_2)),
\end{equation}
while the right hand side can be written as
\begin{equation}
\overline{\mu}'_\Lambda(\lambda_1\otimes f(a_1),\lambda_2\otimes f(a_2))=%
\lambda_2\lambda_1\otimes\mu'(f(a_1),f(a_2)),
\end{equation}
where $f,\mu,\mu'$ are the maps of $\fieldK$-super vector spaces corresponding to the barred versions by
Prop.~\ref{isokmod}. Thus these maps are $\fieldK$-super algebra homomorphisms.
\end{proof}

This argument applies to every
$\fieldK$-multilinear superalgebraic structure which is defined by relations involving finitely 
many arguments, e.g., associative superalgebras \cite{M:Infinite}.

\begin{dfn}
A $\olk$-module $\cv$ in $\catsets^\catgr$ is called superrepresentable if it is isomorphic to
$\olv$ for some $\fieldK$-super vector space $V$.
\end{dfn}

Due to Prop.~\ref{isokmod}, the superrepresentable $\olk$-modules form a full subcategory in 
$\catmod_{\olk}(\catsets^\catgr)$, and $\olcdot$ is an equivalence between this subcategory and
$\catsmod_\fieldK(\catsets)$. Non-superrepresentable $\olk$-modules do indeed exist. An example that
will prove useful later on is $\olv^{nil}$. Let $V$ be some $\fieldK$-super vector space
and let $\Lambda^{nil}$ be the nilpotent ideal in $\Lambda$. Then one can define a $\olk$-module by
setting
\begin{eqnarray}
\label{vnildef}
\olv^{nil}(\Lambda) &:=& (\Lambda^{nil}\otimes_\fieldK V)_{\bar{0}}\\
\olv(\varphi) &:=& (\varphi\otimes\mathrm{id}_{V})\big|_{\olv^{nil}(\Lambda)}\qquad\mathrm{for}\quad%
\varphi:\Lambda\to\Lambda'.
\end{eqnarray}
For every $\Lambda$, one has
\begin{equation}
\label{eq:nildec}
\olv(\Lambda)=V_{\bar{0}}\oplus\olv^{nil}(\Lambda)=V_{\bar{0}}\oplus%
\left(\Lambda_{\bar{0}}^{nil}\otimes V_{\bar{0}}\right)\oplus%
\left(\Lambda_{\bar{1}}^{nil}\otimes V_{\bar{1}}\right).
\end{equation}
Hence $\olv^{nil}$ is superrepresentable if and only if $V_{\bar{0}}=0$, in which case $V$ itself
superrepresents it.

Finally, the change of parity functor $\Pi$ can be carried over to the category of 
superrepresentable $\olk$-modules in an obvious way.

\begin{dfn}
\label{def:copf2}
Let $\mathsf{SRepMod}_{\olk}\subset\catsmod_{\olk}$ be the full subcategory of superrepresentable 
$\olk$-modules. The change of parity functor $\olpi$ is defined as
\begin{eqnarray}
\olpi:\mathsf{SRepMod}_{\olk} &\to& \mathsf{SRepMod}_{\olk}\\
\olv &\mapsto& \overline{(\Pi(V))}.
\end{eqnarray}
\end{dfn}

Prop.~\ref{isokmod} also ensures the existence of inner Hom objects in
the category $\mathsf{SRepMod}_{\olk}$. As the following Corollary shows, we may choose these 
to be just the ``barred'' versions of the inner Hom objects in $\catsmod_\fieldK(\catsets)$.

\begin{cor}
\label{cor:ihomsrep}
Let $\olv_1,\ldots,\olv_n,\olv$ be superrepresentable $\olk$-modules. Then there exists a
superrepresentable inner Hom object
\begin{equation}
\cl^n_{\olk}(\olv_1,\ldots,\olv_n;\olv)\cong\overline{\ihom(V_1,\ldots,V_n;V)}=\overline{L^n_{\fieldK}(V_1,\ldots,V_n;V)}.
\end{equation}
\end{cor}
\begin{proof}
If a superrepresentable inner Hom object exists, it has to satisfy
\[
L^1_{\olk}(\olw;\cl^n_{\olk}(\olv_1,\ldots,\olv_n;\olv))\cong L^{n+1}_{\olk}(\olw,\olv_1,\ldots,\olv_n;\olv)
\]
for all superrepresentalble $\olk$-modules $\olw$. By Prop.~\ref{isokmod}, the space on the right hand
side is isomorphic to
\[
L^{n+1}_{\fieldK,\bar{0}}(W,V_1,\ldots,V_n;V)\cong L^1_{\fieldK,\bar{0}}(W,\ihom(V_1,\ldots,V_n;V))
\]
because inner Hom objects exist in $\catsvect_{\fieldK}$.
This in turn is, by Prop.~\ref{isokmod}, isomorphic to
\[
L^1_{\olk}(\olw,\overline{\ihom(V_1,\ldots,V_n;V)}),
\]
so $\overline{\ihom(V_1,\ldots,V_n;V)}=\overline{L^n_{\fieldK}(V_1,\ldots,V_n;V)}$ satisfies the conditions of an inner Hom object in 
$\mathsf{SRepMod}_{\olk}$.
\end{proof}

This means we may choose the internal Hom-functors
\[
\cl^n_{\olk}:(\catmod_{\olk}^\circ)^n\times\catmod_{\olk}\to\catmod_{\olk}
\]
to be given by the prescription
\[
\cl^n_{\olk}(\olv_1,\ldots,\olv_n;\olv)=\overline{L^n_{\fieldK}(V_1,\ldots,V_n;V)}
\]
for all super vector spaces $V_1,\ldots,V_n,V$, and we will do so in the following.

In a completely analogous way one shows that $\mathsf{SRepMod}_{\olk}$ possesses tensor products over
$\olk$, and that these tensor products inherit all the usual properties from $\catsvect_\fieldK$.

\subsection{Generators for $\catsets^\catgr$ and a criterion for superrepresentability}

Before displaying some examples for the use of the categorical formulation we will construct a
set of generators for $\catsets^\catgr$. This set will turn out to generate many of the 
subcategories of $\catsets^\catgr$ we are
interested in, and moreover will provide an interpretation of the point sets produced by the functor of
points used in supergeometry.

Consider the functors $\cp(\Lambda)$ in $\catsets^\catgr$ defined by
\begin{eqnarray}
\cp(\Lambda):\Lambda' &\mapsto& \Hom_{\catsalg_{\fieldK}}(\Lambda,\Lambda')\\
\nonumber
(\varphi:\Lambda'\to\Lambda'') &\mapsto& \left(\begin{array}{ccc}
\cp(\Lambda)(\varphi):\cp(\Lambda)(\Lambda') &\to& \cp(\Lambda)(\Lambda'')\\
u &\mapsto& \varphi\circ u\end{array}\right).
\end{eqnarray}
From Prop.~\ref{homgr}, we see that $\cp(\Lambda_n)\cong\ovlr{0|n}$.
The functors $\cp(\Lambda)$ are thus just the superpoints defined in Section \ref{sect:spoint}.

\begin{lemma}
\label{sgrgen}
For any functor $F\in\catsets^\catgr$ one has a bijection $F(\Lambda)\cong\Hom(\cp(\Lambda),F)$.
\end{lemma}
\begin{proof}
This is just the statement of the Yoneda lemma for the catgory $\catgr$, which asserts that the map
\begin{eqnarray}
\Hom_{\catsets^\catgr}(\Hom(\Lambda,-),F) &\to& F(\Lambda)\\
\eta &\mapsto& \eta_\Lambda(\id_\Lambda),
\end{eqnarray}
where $\eta_\Lambda$ is the $\Lambda$-component of the functor morphism $\eta$, is a bijection.
\end{proof}

We recall the following definition.

\begin{dfn}
Let $\catc$ be a category. A set $\{G_i\in\catc\}_{i\in I}$ of objects for some index set $I$ is called a set of generators
for $\catc$ if for any pair $f,g:A\to B$ of distinct morphisms between $A,B\in\catc$ 
there exists $i\in I$ and a morphism $s:G_i\to A$ such that the compositions
\begin{equation}
\xymatrix{G_i \ar[r]^s & A \ar@<0.5ex>[r]^f \ar@<-0.5ex>[r]_g & B}
\end{equation}
are still distinct.
\end{dfn}

A set of generators is thus able to keep all distinct morphisms distinct under precomposition. Clearly, given a functor
$F:\catc\to\catsets$, it is enough to consider all the $G_i$-points of it in order to find out whether it
is representable. Obviously, if a generator set exists, then it need not at all be unique.

\begin{cor}
The set $\{\cp(\Lambda)\mid \Lambda\in\catgr\}$ is a set of generators for $\catsets^\catgr$.
\end{cor}
\begin{proof}
Let $\Phi,\Psi:F\to G$ be two distinct morphisms in $\catsets^\catgr$, i.e., there exists a 
$\Lambda\in\catgr$ and $x\in F(\Lambda)$ such that $\Phi_\Lambda(x)\neq\Psi_\Lambda(x)$. Then any
morphism $\eta:\cp(\Lambda)\to F$ which maps $\eta_\Lambda(\id_\Lambda)=x$ separates $\Phi$ and $\Psi$.
\end{proof}

Note, however, that although the $\cp(\Lambda)$ are superrepresentable $\olk$-mod\-ules, they do not
generate the subcategory $\catmod_{\olk}(\catsets^\catgr)$ \emph{as $\olk$-modules} (of course, they
generate it as a restriction of $\catsets^\catgr$). The reason is that the additional structure of a
$\olk$-module restricts the morphisms in such a way that the superpoints are not sufficient anymore to
reproduce the $\Lambda$-points. For example, for any super vector space $V$ over $\fieldK$ we have
\[
\Hom_{\catmod_{\olk}(\catsets^\catgr)}(\cp(\olk),\olv)\cong\Hom_{\catsvect_\fieldK}(\fieldK^{0|0},V)=\{0\},
\]
by Prop.~\ref{isokmod} while
\[
\Hom_{\catsets^\catgr}(\cp(\olk),\olv)\cong V_{\bar{0}}
\]
as sets.

We can obtain a statement analogous to Lemma \ref{sgrgen} by
considering the Grassmann superalgebras
as super vector spaces and forming their functors $\overline{\Lambda}$. Then we obviously have
\begin{equation}
\Hom(\overline{\Lambda},\olv)\cong\Hom_{\catsvect_\fieldK}(\Lambda,V)\cong%
(\Lambda_{\bar{0}}\otimes V_{\bar{0}})\oplus(\Lambda_{\bar{1}}\otimes V_{\bar{1}})\cong\olv(\Lambda).
\end{equation}
The following Lemma shows that there is a simpler possibility.
\begin{lemma}
\label{gensvect}
$\{\fieldK^{1|1}\}$ is a generator set for $\catsvect_\fieldK$.
\end{lemma}
\begin{proof}
It is well known that $\{\fieldK^{1}\}$ is a generator set for the category of
$\fieldK$-vector spaces.
Two morphisms of super vector spaces may differ either in an even or an odd subspace, so
one even and one odd dimension is sufficient to separate all distinct morphisms between super vector spaces.
\end{proof}

Prop.~\ref{isokmod} then implies that $\overline{\fieldK^{1|1}}$ is sufficient to generate 
$\catsrepmod_{\olk}(\catsets^\catgr)$, by
\begin{equation}
\Hom_{\catsrepmod_{\olk}}(\overline{\fieldK^{1|1}},\olv)\cong\Hom_{\catsvect_\fieldK}(\fieldK^{1|1},V)\cong
V_{\bar{0}}\oplus V_{\bar{1}}.
\end{equation}

The fact that $\overline{\Lambda}_1\cong\ovlr{1|1}$ is already a set of generators shows that
superrepresentability is a severe restriction on a $\olk$-module.
\begin{prop}
A $\olk$-module $\cv\in\catsets^\catgr$ is superrepresentable if and only if
\begin{enumerate}
\item the terminal morphism $\cv(\epsilon_{\Lambda_1}):\cv(\Lambda_1)\to\cv(\fieldK)$ is surjective and
\item one has
\begin{equation}
\cv\cong\overline{\cv(\fieldK)\oplus\Pi(\ker\,\cv(\epsilon_{\Lambda_1}))}.
\end{equation}
\end{enumerate}
\end{prop}
\begin{proof}
That follows directly from the definition a superrepresentable module and from (\ref{def:olv}).
\end{proof}

For non-superrepresentable modules, $\overline{\fieldK^{1|1}}$ is in general not a generator set.

\subsection{Application: endomorphisms and supertraces}
\label{sect:endotr}

As a small demonstration of the uses of the categorical approach, we will show how the definition of
the supertrace of an endomorphism of a $\fieldK$-super vector space $V$ is induced from the braiding
isomorphisms in $\catsvect_\fieldK$.

The endomorphisms $\End(V)=\Hom(V,V)$ are, as remarked above, only a $\fieldK$-vector space. In the
spirit of the functor of points, we should rather study the inner Hom-object 
$\iend(V)=\ihom(V,V)$. Suppose we were not able to determine this super vector space by direct
manipulations. Then we could construct it using $\{\overline{\Lambda}\}$ as a set of generators:
\[
\ihom(\olv,\olv)(\Lambda)=\Hom(\overline{\Lambda}\otimes \olv,\olv)\cong%
(\Lambda_{\bar{0}}\otimes\Hom(V,V))\oplus(\Lambda_{\bar{1}}\otimes\Hom(\Pi(V),V)).
\]
These are clearly the points of a superrepresentable module associated to
\[
\iend(V)=\Hom(V,V)\oplus\Pi(\Hom(\Pi(V),V)),
\]
which coincides, of course, with the common definition as the superspace of even and odd 
linear morphisms $V\to V$. Its elements may be interpreted as matrices whose diagonal 
$V_{\bar{0}}\times V_{\bar{0}}$- and $V_{\bar{1}}\times V_{\bar{1}}$-blocks contain even entries
(i.e., are from $\fieldK$) and
the off-diagonal blocks contain odd entries (are from $\Pi(\fieldK)$). The sets
$\ihom(\olv,\olv)(\Lambda)$ can be interpreted in the very same way, with the even blocks then containing
elements of $\Lambda_{\bar{0}}$ and the odd blocks elements of $\Lambda_{\bar{1}}$. 

The dual super vector space $V^*$ is by definition
\[
V^*:=\ihom(V,\fieldK),
\]
so it is
\begin{eqnarray*}
\olv^*(\fieldK) &\cong& \Hom_{\catsvect_{\fieldK}}(V,\fieldK)=V_{\bar{0}}^*\\
\olv^*(\Lambda_1) &\cong& \Hom_{\catsvect_{\fieldK}}(V\oplus\Pi(V),\fieldK)=%
V_{\bar{0}}^*\oplus V_{\bar{1}}^*,
\end{eqnarray*}
where $V_{\bar{1}}^*$ is the ordinary dual space, i.e., 
$V_{\bar{1}}^*=\Hom_{\catvect_\fieldK}(V_{\bar{1}},\fieldK)$.
One easily checks that the higher points are just those of a superrepresentable module defined by
\[
V^*=V_{\bar{0}}^*\oplus\Pi(V_{\bar{1}}^*).
\]
Now on one hand, we have an isomorphism
\begin{eqnarray}
\phi:V\otimes V^* &\to& \iend(V)\\
\nonumber
v\otimes f &\mapsto& (w\mapsto f(w)\cdot v).
\end{eqnarray}
On the other hand, we have the evaluation map
\begin{eqnarray}
\mathrm{ev}:V^*\otimes V &\to& \fieldK\\
f \otimes v &\mapsto& f(v).
\end{eqnarray}
The trace is generally defined to be the composition
\begin{equation}
\Tr:\iend(V)\stackrel{\phi^{-1}}{\longrightarrow}V\otimes V^*
\stackrel{c_{V,V^*}}{\longrightarrow}V^*\otimes V\stackrel{\mathrm{ev}}{\longrightarrow}\fieldK
\end{equation}
Now the minus sign occuring in the super trace is induced from the braiding isomorphism $c_{V,V^*}$.
Let $\{e_i\}$ be a basis of $V$ and $\{e^j\}$ be the
dual basis of $V^*$. Then for $\psi\in\iend(V)$, $\phi^{-1}(\psi)$ is simply its decomposition into the basis,
$\psi^i_j e_i\otimes e^j$. Applying $c_{V,V^*}$ yields $(-1)^{p(e_i)p(e^j)}\psi^i_j e^j\otimes e_i$ and thus,
$\Tr(\psi)=(-1)^{p(e_i)}\psi^i_i$.

This example shows that the
categorical approach can be a guideline how to ``superize'' notions from algebra if guessing them
is impossible. It is also clear that the direct use of odd and even elements can be more elegant and
transparent, if one deals with a finite-dimensional or linear problem. But especially in geometry, 
as will be discussed next, the categorical method can sometimes be the only viable one.

\section{Supergeometry in the categorical setting}

We will now proceed to define superdomains and supermanifolds using the categorical formalism
introduced in the previous section.
Superdomains will be defined as open subfunctors of
superrepresentable $\olk$-modules. Since we can define $\fieldK$-super vector spaces and thus $\olk$-modules
of arbitrary dimensions, infinite-dimensional superdomains will also become available. These open
subfunctors will be glued in a straightforward way to form
(possibly infinite-dimensional) supermanifolds.

To make the notion of an open subobject sensible, one has to define an analogue of a topology on 
$\olk$-modules. Since the category
$\catmod_{\olk}(\catsets^\catgr)$ is not concrete, we will have to resort to the more abstract concept of
a Grothendieck topology (for background, see e.g., \cite{FGIKNV:Fundamental}, \cite{SGA4}, \cite{MM:Sheaves}). 

\subsection{The topology on $\cattop^\catgr$}

To define supermanifolds, the only functors of $\catsets^\catgr$ which are of relevance
are superrepresentable $\olk$-modules. But these are not merely functors into $\catsets$, 
but actually into the
category $\cattop$ of topological spaces if we assume that the even and odd parts of the super vector spaces
$V$ that are used can be given the structure of topological vector spaces.

\begin{dfn}
Let $\cf,\cf'$ be functors in $\cattop^\catgr$. $\cf'$ is called a subfunctor of $\cf$, if
\begin{enumerate} 
\item for every $\Lambda\in\catgr$, $\cf'(\Lambda)$ is a topological subspace of $\cf(\Lambda)$, and
\item the family of inclusions $\{\cf'(\Lambda)\subset\cf(\Lambda)\big|\Lambda\in\catgr\}$ forms a
functor morphism.
\end{enumerate}
In this case, one just writes $\cf'\subset\cf$. $\cf'$ is called an open subfunctor of $\cf$ if, in 
addition, each $\cf'(\Lambda)$ is open in $\cf(\Lambda)$.
\end{dfn}

\begin{dfn}
Let $\cf',\cf''$ be subfunctors of $\cf\in\cattop^\catgr$. Then the intersection $\cf'\cap\cf''$ is
the functor whose points are
\begin{equation}
(\cf'\cap\cf'')(\Lambda):=\cf'(\Lambda)\cap\cf''(\Lambda).
\end{equation}
The union $\cf'\cup\cf''$ is the functor defined by
\begin{equation}
(\cf'\cup\cf'')(\Lambda):=\cf'(\Lambda)\cup\cf''(\Lambda).
\end{equation}
A morphism $\varphi:\Lambda\to\Lambda'$ is mapped by
$\cf'\cap\cf''$ resp. $\cf'\cup\cf''$ to the corresponding restrictions 
$\cf(\varphi)\big|_{\cf'(\Lambda)\cap\cf''(\Lambda)}$ resp.
$\cf(\varphi)\big|_{\cf'(\Lambda)\cup\cf''(\Lambda)}$.
\end{dfn}

Clearly, both $\cf'\cap\cf''$ and $\cf'\cup\cf''$ are again subfunctors of $\cf$. The functor
$\mathrm{emp}:\Lambda\to\emptyset$ is the initial object in $\cattop^\catgr$.

\begin{dfn}
A functor morphism $g:\cf''\to\cf$ is called open if there exists a factorization
\[
\begin{CD}
g:\cf'' @>f>> \cf'\subset\cf
\end{CD}
\]
such that $f$ is an isomorphism and $\cf'$ is an open subfunctor of $\cf$.
\end{dfn}

Finally, the notion of an open covering can be carried over straightforwardly.

\begin{dfn}
A family
$\{u_\alpha:\cu_\alpha\to\cf\}$ of open functor morphisms is called an open covering of $\cf$ if for
each $\Lambda\in\catgr$, the family of maps
\[
u_{\alpha,\Lambda}:\cu_\alpha(\Lambda)\to\cf(\Lambda)
\]
is an open covering of the topological space $\cf(\Lambda)$.
\end{dfn}

It is obvious from the definitions that this definition of open coverings
endows $\cattop^\catgr$ with a Grothendieck topology, namely by simply pulling back the global classical
topology from $\cattop$. We will from now on always assume that $\cattop^\catgr$ is endowed with this
particular topology.

An example of an open subfunctor can be constructed in the following way: let  $\cf\in\cattop^\catgr$ be
an arbitrary functor, and let $U\subset\cf(\realR)$ be an open subset of its underlying set. Then we can
construct an open subfunctor $\cu\subset\cf$ by setting
\begin{eqnarray}
\nonumber
\cu(\realR) &:=& U\\
\label{restfunc}
\cu(\Lambda) &:=& \cf(\epsilon_\Lambda)^{-1}(U)\subset\cf(\Lambda)\\
\nonumber
\cu(\varphi) &:=& \cf(\varphi)\big|_{\cu(\Lambda)}\qquad\mathrm{for}\,\,\varphi:\Lambda\to\Lambda'.
\end{eqnarray}
Here, $\epsilon_\Lambda:\Lambda\to\realR$ is the terminal morphism of $\Lambda\in\catgr$.
It is clear from the definition that the inclusion $\cu\subset\cf$ is indeed a functor morphism. We will
denote subfunctors of this form by $\cu=\cf\big|_U$ and call them restrictions. It will turn out that only
such subfunctors qualify as ``affine'' domains for the construction of supermanifolds.

\subsection{Superdomains and supersmooth morphisms}

Now everything is prepared for the introduction of (possibly in\-fin\-ite-di\-men\-sion\-al) superdomains.

\begin{dfn}
Let $\cv$ be a superrepresentable $\olk$-module in $\cattop^\catgr$. $\cv$ will be called a locally
convex, resp.~Fr\'echet, resp.~Banach $\olk$-module if for every $\Lambda\in\catgr$, the topological
vector space $\cv(\Lambda)$ is locally convex, resp.~Fr\'echet, resp.~Banach.
\end{dfn}

\begin{dfn}
An open subfunctor $\cf$ of a locally convex (resp.~Fr\'echet, resp.~Banach) $\olk$-module in
$\cattop^\catgr$ is called a real (or complex, whichever $\olk$ is) locally convex 
(resp.~Fr\'echet, resp.~Banach) superdomain.
\end{dfn}

\begin{dfn}
A functor $\cf\in\cattop^\catgr$ is called locally isomorphic to real (or complex) locally convex
(resp.~Fr\'echet, resp.~Banach) superdomains if there exists an open covering 
$\{u_\alpha:\cu_\alpha\to\cf\}$ of $\cf$ such that each $\cu_\alpha$ is a locally convex
(resp.~Fr\'echet, resp.~Banach) superdomain.
\end{dfn}

Restrictions, as it turns out, are the only open subfunctors a superrepresentable $\olk$-module has.

\begin{prop}
\label{subf}
Any open subfunctor $\cu\subset\olv$ of a Banach (resp.~Fr\'echet, resp.~locally convex) $\olk$-module $\olv$
is a restriction $\olv\big|_U$, where $U=\cu(\fieldK)$ (cf. (\ref{restfunc})).
\end{prop}
\begin{proof}
Clearly, one can write 
\begin{equation}
\olv=\olv\big|_{V_{\bar{0}}}
\end{equation}
for a superrepresentable $\olk$-module $\olv$ which is represented by $V$, since
\begin{equation}
\olv(\Lambda)=\olv(\epsilon_\Lambda)^{-1}(V_{\bar{0}}).
\end{equation}
Let now $\cu$ be an arbitrary open subfunctor of $\olv$ and $U\subset V_{\bar{0}}$ be its $\fieldK$-points.
The inclusion $\cu\subset\olv$ must be a functor morphism, therefore the diagram
\begin{equation}
\begin{CD}
\cu(\Lambda) @>\subset>> \olv(\Lambda)\\
@VV{\cu(\epsilon_\Lambda)}V @VV{\olv(\epsilon_\Lambda)}V\\
U @>\subset>> \olv(\fieldK)=V_{\bar{0}}
\end{CD}
\end{equation}
has to commute for all $\Lambda\in\catgr$. This enforces
\begin{equation}
\cu(\Lambda)=\olv(\epsilon_\Lambda)^{-1}(U)
\end{equation}
for all $\Lambda$. For any morphism $\varphi:\Lambda\to\Lambda'$ of Grassmann algebras, we also have
\begin{equation}
\begin{CD}
\cu(\Lambda) @>\subset>> \olv(\Lambda)\\
@VV\cu(\varphi)V @VV\olv(\varphi)V\\
\cu(\Lambda') @>\subset>> \olv(\Lambda')
\end{CD},
\end{equation}
which commutes again, because the inclusion is a functor morphism. So,
\begin{equation}
\cu(\varphi)=\olv(\varphi)\big|_{\cu(\Lambda)}.
\end{equation}
\end{proof}

By the properties of unions and intersections of open subfunctors, we obtain the following Corollary.

\begin{cor}
Let $\cf\in\cattop^\catgr$ be a functor which is locally isomorphic to Banach 
(resp.~Fr\'echet, resp.~locally convex) superdomains. Then every open subfunctor $\cu\subset\cf$ is 
a restriction $\cf\big|_U$.
\end{cor}
\begin{proof}
Let $\{u_\alpha:\cu_\alpha\to\cf\}$ be an open cover of $\cf$ by superdomains of the appropriate type
and let $\cu$ be an arbitrary open subfunctor of $\cf$. Then
every intersection $\cu\cap\cu_\alpha$ is a superdomain, i.e. is a restriction $\cf\big|_{U\cap U_\alpha}$,
where $U,U_\alpha$ are the underlying open sets of the respective functors. Therefore,
\[
\cu(\Lambda)=\bigcup_{\alpha}(\cf(\epsilon_\Lambda))^{-1}(U_\alpha\cap U)=(\cf(\epsilon_\Lambda))^{-1}(U)
\]
By the same argument as in Prop.~\ref{subf} (functoriality of inclusions),
the images of the morphisms of $\catgr$ under $\cu$ must be the restrictions of those of 
$\cf$.
\end{proof}

For the rest of this section, we will focus exclusively on
Banach superdomains because Banach spaces are analytically nicest. After the necessary modifications,
most of the constructions described below will carry over at least to the tame Fr\'echet case.

\begin{dfn}
\label{bsdom}
Let $\cv\big|_U$ and $\cv'\big|_{U'}$ be two real Banach superdomains. A functor morphism 
$f:\cv\big|_U\to\cv'\big|_{U'}$ is called supersmooth if 
\begin{enumerate}
\item the map
\begin{equation}
f_\Lambda:\cv\big|_U(\Lambda)\to\cv'\big|_{U'}(\Lambda)
\end{equation}
is smooth for every $\Lambda\in\catgr$, and
\item for every $u\in \cv\big|_U(\Lambda)$, the derivative
\begin{equation}
Df_\Lambda(u):\cv(\Lambda)\to\cv'(\Lambda)
\end{equation}
is $\Lambda_{\bar{0}}$-linear.
\end{enumerate}
\end{dfn}

The second condition is necessary and sufficient to turn the sets of differential morphisms
\begin{eqnarray}
(Df)_\Lambda:\cv\big|_U(\Lambda)\times\cv(\Lambda) &\to& \cv'(\Lambda)\\
(Df)_\Lambda(u,v) &=& (Df_\Lambda(u))(v)
\end{eqnarray}
into a $\cv\big|_U$-family of $\olr$-linear morphisms $\cv\to\cv'$. This, in turn, has to be required to
make the differentiable structures on the sets $\cv\big|_U(\Lambda)$, $\cv'\big|_{U'}(\Lambda)$
functorial with respect to $\Lambda$, i.e., compatible with morphisms $\varphi:\Lambda\to\Lambda'$.

Together with supersmooth morphisms, smooth Banach superdomains form a category $\catbsdom$. Replacing 
``smooth'' with ``real analytic'' in 
Def.~\ref{bsdom} leads to the definition of the category of real analytic
superdomains. 

For complex analytic domains, there seem to be two different approaches at first. One can start with superdomains
in $\cattop^\catgr$ which are isomorphic to open subfunctors of superrepresentable $\olc$-modules, and
define morphisms to be complex superanalytic functor morphisms. On the other hand, one could as well use
the category $\catgr^\complexC$ of Grassmann algebras over $\complexC$ from the very beginning on, studying only functors in 
$\cattop^{\catgr^\complexC}$ and using morphisms which are analytic in their complex coordinates. 
However, the two resulting categories are equivalent: the $\Lambda$-points of some superrepresentable
$\olc$ module $V$ are $(V\otimes_\realR\Lambda)_{\bar{0}}$. But
\begin{equation}
(V\otimes_\realR\Lambda)_{\bar{0}}\cong%
(V^\realR\otimes_\realR\complexC\otimes_\realR\Lambda)_{\bar{0}}\cong%
(V^\realR\otimes_\realR\Lambda^\complexC)_{\bar{0}},
\end{equation}
where $V^\realR$ is the real super vector space underlying $V$. The notion of supersmoothness will be
identical for both cases: if we use real Grassmann algebras but $\olc$-modules, we will require the component
maps of a functor morphism between superdomains to be holomorphic, which will imply that the differentials are
$\olc$-linear. If we use complex Grassmann algebras, $\Lambda_{\bar{0}}$-linearity of the differentials will
imply $\complexC$-linearity and thus holomorphicity of the component maps of morphisms between superdomains.
We will use the term ``supersmooth'' in this sense from now on, i.e., denoting holomorphicity in the complex
case.

\subsection{The fine structure of supersmooth morphisms}

An explicit description of supersmooth morphisms and a means to determine the sets 
$\Hom_\catbsdom(\cu,\cv)$ of
morphisms between Banach superdomains will be of great value for the following considerations.

\begin{thm}[see also \cite{M:Infinite}]
\label{thm:skel}
A morphism $f:\olv\big|_U\to\olv'\big|_{U'}$ of Banach superdomains is supersmooth if and only if
there exists a smooth map
\[
f_0:U\to U'
\]
and for all $k\geq 1$ smooth maps
\[
f_k:U\to \Sym^k(V_{\bar{1}};V')
\]
such that for every component $f_\Lambda$ of $f$, one has
\begin{equation}
\label{skeleton}
f_\Lambda(u+v_0^{nil}+v_1^{nil})=\sum_{k,m=0}\frac{1}{k!m!}%
\overline{D^kf_m(u)}_\Lambda(\underbrace{v_0^{nil},\ldots,v_0^{nil}}_{k\,\,times},%
\underbrace{v_1^{nil},\ldots,v_1^{nil}}_{m\,\,times}).
\end{equation}
The collection of maps $\{f_k\mid k\geq 0\}$ will be denoted as $f_\bullet$ and will be called a 
skeleton of the morphism $f$. A skeleton, if it exists, is unique.
\end{thm}

Some explanations are in order.
The argument $v\in V(\Lambda)$ of $f_\Lambda$ has been decomposed according to (\ref{eq:nildec}) into
its underlying part $u\in V_{\bar{0}}$ and the parts $v_0^{nil}$ and $v_1^{nil}$, which
contain even and odd elements of $\Lambda$, respectively. This decomposition is unique,
cf. (\ref{eq:nildec}). The k-th differential $D^kf_m(u)$ is, for each $u\in U$, a symmetric element of
$L^{k+m}_{\bar{0}}(V_{\bar{0}}^{\otimes k},V_{\bar{1}}^{\otimes m};V')$, i.e., an even linear map which is supersymmetric within its
first $k$ and within its last $m$ arguments.
Then $\overline{D^kf_m(u)}_\Lambda$ is the $\Lambda$-component
of the associated $\olk$-multilinear functor morphism. This functor morphism is uniquely determined
by Prop.~\ref{isokmod}.

\begin{proof}
Let $f:\olv\big|_U\to\olv'\big|_{U'}$ be a given supersmooth morphism between two Banach superdomains.
We have to show that $f$ gives rise to a unique skeleton $f_\bullet$.

Clearly, $f_\fieldK$ is a smooth map $U\to U'$ and therefore qualifies as the required map $f_0$.
Let $f_\Lambda$ be one of the higher components of $f$. By assumption, it is a smooth map between the Banach
domains $\olv\big|_U(\Lambda)$ and $\olv'\big|_{U'}(\Lambda)$. Let $v=u+v_0^{nil}+v_1^{nil}$ be a point of 
$\olv\big|_U(\Lambda)$.
We can use the Taylor expansion around $u$ to write
\begin{equation}
\label{taylor}
f_\Lambda(u+v_0^{nil}+v_1^{nil})=\sum_{k=0}\frac{1}{k!}%
D^kf_\Lambda(u)(\underbrace{v_0^{nil}+v_1^{nil},\ldots,v_0^{nil}+v_1^{nil}}_{k\,\,times}),
\end{equation}
where the differentials $D^kf_\Lambda$ are $U$-families of symmetric $\fieldK$-linear maps $V(\Lambda)^k\to V'(\Lambda)$. 
Using linearity and symmetry, we can rewrite this as
\begin{equation}
\label{taylor2}
f_\Lambda(u+v_0^{nil}+v_1^{nil})=\sum_{k=0,m\leq k}\frac{1}{m!(k-m)!}%
D^kf_{\Lambda}(u)(\underbrace{v_0^{nil},\ldots,v_0^{nil}}_{(k-m)\,\,times},%
\underbrace{v_1^{nil},\ldots,v_1^{nil}}_{m\,\,times}).
\end{equation}
Collecting all terms which depend on m copies of $v_1^{nil}$ produces the sum
\begin{equation}
\label{preon}
\sum_{k=0}\frac{1}{k!m!}D^{k+m}f_\Lambda(u)(\underbrace{v_0^{nil},\ldots,v_0^{nil}}_{k\,\,times},%
\underbrace{v_1^{nil},\ldots,v_1^{nil}}_{m\,\,times}).
\end{equation}
By definition, the differentials of the component maps of a supersmooth morphism are
$\Lambda_{\bar{0}}$-linear. Therefore, the maps
\begin{equation}
\label{deffm}
f_{\Lambda,m}(u) := D^mf_\Lambda(u)\big|_{\olv_{\bar{1}}(\Lambda)}:U\to \Sym^m(\olv_{\bar{1}}(\Lambda);\olv'(\Lambda))
\end{equation}
form the $\Lambda$-components of a smooth $U$-family of $\olk$-linear morphisms $\olf_m(u):\olv_{\bar{1}}^{\otimes m}\to\olv'$. By
Prop. (\ref{isokmod}), $\olf_m$ corresponds to a unique family of $\fieldK$-linear supersymmetric maps
$f_m:V_{\bar{1}}^{\otimes k}\to V'$. Then $D^kf_m$ is, in turn, a smooth $U$-family of maps
$V_{\bar{0}}^{\otimes k}\times V_{\bar{1}}^{\otimes m}\to V'$, which are supersymmetric within the two groups
of arguments, and to which we assign a unique $U$-family $\overline{D^kf_m}$. By uniqueness, the $\Lambda$-components of
$\overline{D^kf_m}$ are identical to the $D^kf_\Lambda$ in eq. (\ref{taylor2}), which also implies that the sums in
eqs. (\ref{taylor}) and (\ref{taylor2}) are finite.

Clearly, the family $\{f_\Lambda\}$ uniquely determines the family $\{f_m\}$, which is a skeleton for $f$.
Conversely, every collection of smooth maps $f_n:U\to\Sym^n(V_{\bar{1}};V')$
obviously is the skeleton of a supersmooth functor morphism $f:\olv\big|_U\to\olv'$, whose components are
defined by (\ref{skeleton}). One also verifies that compositions of skeletons are skeletons for the composed
maps.
\end{proof}

As pointed out in \cite{M:Infinite}, one could use this result to define Banach superdomains
completely without the use of functors and Grassmann algebras, namely as pairs 
$(U\subset V_{\bar{0}},V_{\bar{1}})$, where both $V_{\bar{0}},V_{\bar{1}}$ are ordinary 
Banach spaces and $U$ is an open domain in $V_{\bar{0}}$, and the morphisms between these pairs are 
skeletons. This is a point of view that we will not pursue here, however.

\subsection{Banach supermanifolds}
\label{sect:bansman}

A superrepresentable $\olk$-module $\olv$ all of whose points $\olv(\Lambda)$ have been
endowed with Banach space structures in a functorial manner is not just a functor $\catgr\to\cattop$, 
but actually a functor $\catgr\to\catman$,
where $\catman$ is the category of smooth Banach manifolds. In order to study Banach supermanifolds,
we will now restrict our attention to the category $\catman^\catgr$.

\begin{dfn}[see also \cite{M:Infinite}]
\label{def:sman}
Let $\cf$ be a functor in $\catman^\catgr$. An open covering $\ca=\{u_\alpha:\cu_\alpha\to\cf\}_{\alpha\in I}$
of $\cf$ is called a supersmooth atlas on $\cf$ if
\begin{enumerate}
\item every $\cu_\alpha$ is a Banach superdomain,
\item for every pair $\alpha,\beta\in I$, the fiber product 
\begin{equation}
\cu_{\alpha\beta}=\cu_\alpha\times_\cf\cu_\beta\in\catman^\catgr
\end{equation}
can be given the structure of a Banach superdomain such that the projections 
$\Pi_\alpha:\cu_{\alpha\beta}\to\cu_\alpha$ and $\Pi_\beta:\cu_{\alpha\beta}\to\cu_\beta$ are supersmooth.
\end{enumerate}
The maps $u_\alpha:\cu_\alpha\to\cf$ are called charts on $\cf$.
\end{dfn}

\begin{dfn}[see also \cite{M:Infinite}]
Two supersmooth atlases $\ca,\ca'$ on the functor $\cf\in\catman^\catgr$ are said to be
equivalent if their union $\ca\cup\ca'$ is again a supersmooth atlas on $\cf$. A supermanifold $\cm$
is a functor in $\catman^\catgr$ endowed with an equivalence class of atlases.
\end{dfn}

The second condition in Definition \ref{def:sman} needs some explanation. One might think at first that
the fiber product of any two supersmooth superdomains is automatically again a superdomain. But we may
only assume here that $\cf$ is a functor in $\catman^\catgr$, and thus we can a priori only 
assume the fiber product to exist in $\catman^\catgr$. The projection morphisms $\Pi_\alpha,\Pi_\beta$
are therefore only guaranteed to be functor morphisms in $\catman^\catgr$. The second condition therefore
requires
the fiber product of $\cu_\alpha,\cu_\beta$ to exist in the subcategory $\catbsdom\subset\catman^\catgr$.
Since there really exist functor isomorphisms in $\catman^\catgr$ which are not supersmooth (compare
Thm.~\ref{thm:skel}), this is not automatic.

\begin{dfn}[see also \cite{M:Infinite}]
\label{dfn:ssmoothmor}
Let $\cm,\cm'$ be Banach supermanifolds. A functor morphism $f:\cm\to\cm'$ is called supersmooth if for
each pair of charts $u:\cu\to\cm$, $u':\cu'\to\cm'$, the pullback 
\begin{equation}
\xymatrix{
\cu\times_{\cm'}\cu' \ar[d]_{\Pi}\ar[rr]^{\Pi'} && \cu' \ar[d]^{u'}\\
\cu \ar[r]^{u} & \cm \ar[r]^{f} & \cm'
}
\end{equation}
can be given the structure of a Banach superdomain such that its projections $\Pi,\Pi'$ are supersmooth.
\end{dfn}

It is clear that the composition of two supersmooth morphisms of Banach supermanifolds is again
supersmooth. Thus Banach supermanifolds form a category $\catbsman$. If nothing else is specified, we
will from now on only write $\catsman$ for the category of Banach supermanifolds, since we restrict ourselves
to them in this work. The set of supersmooth morphisms $f:\cm\to\cm'$ will be denoted as
$SC^\infty(\cm,\cm'):=\Hom_{\catbsman}(\cm,\cm')$.

According to Prop.~\ref{subf}, every open submanifold of a Banach supermanifold $\cm$ is of the form
$\cu=\cm\big|_U$, where $U$ is an open submanifold of the underlying manifold $M=\cm(\fieldK)$. Clearly,
the restriction of the supermanifold structure of $\cm$ to $\cu$ induces on the latter a
unique supermanifold structure which makes the inclusion $\cu\subset\cm$ a supersmooth morphism.

\begin{rmk}
\label{rmk:smanpoint}
The above definition implies that $M=\cm(\realR)$ is a smooth
manifold and every point set $\cm(\Lambda)$ has the structure of a smooth vector bundle over $M$, with projection
given by the image $\cm(\epsilon_{\Lambda})$ of the terminal morphism $\epsilon_{\Lambda}:\Lambda\to\realR$:
\begin{equation}
\cm(\epsilon_\Lambda):\cm(\Lambda)\to\cm(\realR)=M.
\end{equation}
If $V$ is the super vector space on which $\cm$ is modeled, then the fiber of this bundle is isomorphic to 
the kernel of the map 
$\olv(\epsilon_\Lambda):\olv(\Lambda)\to V_{\bar{0}}$, which is just $\olv^{nil}(\Lambda)$ (cf.~eq.~(\ref{vnildef})).
As we will discuss below, the similarity of this picture with a de Witt supermanifold is not an accident.
\end{rmk}

\subsection{The topology on $\catbsdom$ and $\catsman$}

The categories $\catbsdom$ and $\catsman$ inherit a topology from $\cattop^\catgr$ in a quite straightforward
way. We note first that we can pull back the Grothendieck topology from $\cattop^\catgr$ to $\catman^\catgr$ 
along the forgetful functor $N:\catman^\catgr\to\cattop^\catgr$. This is equivalent to saying that a family
$\{u_\alpha:\cu_\alpha\to\cu\}$ of morphisms in $\catman^\catgr$ is an open covering of $\cu$ if 
and only if it is an open covering of $\cu$
when the $u_\alpha$ are considered as morphisms between objects in $\cattop^\catgr$. In exactly the same
way $\catbsdom$ and $\catsman$ are supposed to inherit their topology as a subcategories of $\catman^\catgr$.

Although we will not make use of it in this article, we want to include in this Section the proof of
a useful fact that will become important, e.g., in the construction of super vector bundles over
supermanifolds.

\begin{prop}
The Grothendieck topology on $\catsman$ is subcanonical, i.e., every representable functor on $\catsman$
is a sheaf.
\end{prop}
\begin{proof}
Let $\cm$ be a given supermanifold and let $\{u_\alpha:\cu_\alpha\to\cu\}$ be an open covering in $\catsman$.
The representable functor $\Hom(-,\cm)$ defines a presheaf on $\catsman$
for which we have to ensure that the gluing axiom holds. Let $f_\alpha:\cu_\alpha\to\cm$ be 
supersmooth morphisms such that for every fibered product $\cu_\alpha\times_\cu\cu_\beta$ the condition
\[
\pi_\alpha^*f_\alpha=\pi_\beta^*f_\beta
\]
holds, where $\pi_{\alpha,\beta}:\cu_\alpha\times_\cu\cu_\beta\to\cu_{\alpha,\beta}$ are the canonical
projections. Then we have to show that there exists a unique supersmooth morphism $f:\cu\to\cm$ such that
\begin{equation}
\label{eq:glue}
u_\alpha^*f=f_\alpha\qquad\forall\alpha.
\end{equation}
Given the $f_\alpha$ we can pointwise, i.e. for every $\Lambda$, patch together the smooth maps
$f_{\alpha\Lambda}$ to yield smooth maps $f_\Lambda:\cu(\Lambda)\to\cm(\Lambda)$ which are 
functorial in $\Lambda$, i.e.,
they form a functor morphism $f$ in $\catman^\catgr$ satisfying (\ref{eq:glue}). It
therefore only remains to be shown that this $f$ is supersmooth. We choose a covering 
$\{v_\beta:\cv_\beta\to\cm\}$ of $\cm$. Then for all $\alpha,\beta$ we have commutative diagrams
\begin{equation}
\xymatrix{\cu_\alpha\times_\cm\cv_\beta \ar[rr]^{\pi_\beta} \ar[d]_{\pi_\alpha} && 
\cv_\beta \ar[d]^{v_\beta} \\
\cu_\alpha \ar[dr]_{u_\alpha} \ar[rr]^{f_\alpha} && \cm\\
& \cu \ar[ru]_f &}
\end{equation}
where the upper square commutes and all its arrows are supersmooth because the $f_\alpha$ are supersmooth 
and the lower triangle commutes by construction. This is precisely the condition for $f$ to be supersmooth,
cf. Def.~\ref{dfn:ssmoothmor}.
\end{proof}

\subsection{Superpoints are generators for $\catsman$}

Superpoints were introduced in section \ref{sect:spoint} as linear supermanifolds corresponding to purely odd
super vector spaces. In Prop.~\ref{spointeq}, it was shown that the category $\catspoint$ is dual to
the category $\catgr$, and a duality was chosen, namely
$\cp:\Lambda \mapsto \Spec(\Lambda)=(\{*\},\Lambda)$.
This amounts to the choice $\cp(\Lambda_n)=\ovlk{0|n}$ for the Grassmann algebra
on $n$ generators over $\fieldK$. But a supermanifold is also a functor 
$\catgr\to\catman$. So, $\cp$ can be considered as a bifunctor
\begin{equation}
\cp:\catgr^\circ\times\catgr\to\catman.
\end{equation}

\begin{prop}
\label{prop:isocp}
There exists an isomorphism of bifunctors
\begin{equation}
\cp\cong\Hom_\catgr(-,-).
\end{equation}
\end{prop}
\begin{proof}
The $\Lambda_m$-points of $\cp(\Lambda_n)$ are
\begin{equation}
\cp(\Lambda_n)(\Lambda_m)=(\Lambda_m\otimes\fieldK^{0|n})_{\bar{0}}%
\cong \Lambda_{m,\bar{1}}\otimes\fieldK^n\cong\Hom_\catgr(\Lambda_n,\Lambda_m),
\end{equation}
where the last isomorphism was proved in Prop.~\ref{homgr}. 
\end{proof}

We will now show that the superpoints generate the category $\catsman$, which basically means that they
play the role for supermanifolds which is played by the manifold $\Spec\fieldK=(\{*\},\fieldK)$ for
ordinary manifolds over $\fieldK$.

\begin{prop}
For $\cm\in\catsman$ and every $\Lambda\in\catgr$, one has a bijection
\begin{equation}
\cm(\Lambda)\cong SC^\infty(\cp(\Lambda),\cm)=\Hom_\catsman(\cp(\Lambda),\cm).
\end{equation}
\end{prop}
\begin{proof}
According to Thm.~\ref{thm:skel}, each supersmooth map $f:\cp(\Lambda_n)\to\cm$ can be expressed in
terms of a skeleton $f_\bullet:\{*\}\to\Sym^\bullet(\realR^{0|n};V')$, where $V'$ is the super vector space,
on which $\cm$ is locally modeled. Now $f_0=f_\fieldK:\{*\}\to M=\cm(\fieldK)$ is a map of the one point set into the
underlying manifold $M$ of $\cm$. Therefore we have a bijection
\[
SC^\infty(\cp(\Lambda_n),\cm)\cong M\times\Sym^{k\geq 1}(\realR^{0|n};V').
\]
But
\[
\Sym^k(\realR^{0|n};V')=\wedge^k(\realR^n,V_{\bar{k}})=\Hom(\wedge^k\realR^n,V_{\bar{k}})\cong%
\Lambda_n^{(k)}\otimes_\fieldK V_{\bar{k}},
\]
where $\wedge^k$ are the alternating $\fieldK$-linear maps on $k$ arguments, $\bar{k}$ is to be
understood mod 2 and $\Lambda_n^{(k)}$ denotes the elements of degree $k$ in $\Lambda_n$. Summing up,
we have found
\[
\Sym^{k\geq 1}(\realR^{0|n};V')\cong {\olv'}^{nil}(\Lambda_n)
\]
and therefore
\begin{equation}
SC^\infty(\cp(\Lambda_n),\cm)\cong M\times {\olv'}^{nil}(\Lambda_n).
\end{equation}
As pointed out in Remark \ref{rmk:smanpoint}, this set is in bijection with $\cm(\Lambda)$.
\end{proof}

The bijection constructed here is clearly functorial both with respect to $\cm$ as well as with respect to
$\Lambda$. This ensures that the set $\{\Spec\Lambda_n=\cp(\Lambda_n)\mid n\in\naturalN\}$ is a set of
generators for $\catsman$. This set is not minimal however, since any of its infinite subsets is a generator
set as well\footnote{I thank V. Molotkov for pointing this out to me.}. A minimal set of generators does
not exist for $\catsman$.

\subsection{Application: Lie supergroups}

Supergroups are a topic where the functor of points has been used already from the very beginning on.
In fact, the desire to construct super analogs of Lie groups was one of the ideas that led to the
introduction of supermanifolds \cite{L:Introduction}.
By definition, a Lie supergroup $\cg$ is a group object in the category of supermanifolds. 
In the language developed above this means that we have morphisms
\begin{eqnarray}
\nonumber
m &:& \cg\times\cg\to\cg\\
i &:& \cg\to\cg\\
\nonumber
e &:& \cp(\realR)\to\cg
\end{eqnarray}
such that each set $\cg(\Lambda)$ becomes a group with group law $m_\Lambda$, inversion $i_\Lambda$ and
unit $e_\Lambda$.
So, a supergroup is not at all just a set with some special group structure,
but rather a tower of groups indexed by $\catgr$ and related by the maps which are induced by functoriality 
with respect to $\catgr$.
Luckily, for matrix supergroups the situation simplifies greatly because they are all subgroups of
general linear supergroups, and these have global coordinates.

For any $\fieldK$-super vector space $V$, the superspace $\iend(V)$ is obviously a $\fieldK$-superalgebra
with multiplication the composition of morphisms. Assuming $V$ to be finite-dimensional and to be
given in some fixed format (i.e., with fixed parities of the basis vectors), we can
express the elements of $\iend(V)$ as super matrices (cf.~Sect.~\ref{sect:endotr}) and the multiplication
\[
\mu:\iend(V)\times\iend(V)\to\iend(V)
\]
is just matrix multiplication. The algebra $\iend(V)$ is not supercommutative, of course.
Switching to the associated $\olk$-module $\ihom(\olv,\olv)$ translates $\mu$
into a functor morphism
\[
\overline{\mu}:\ihom(\olv,\olv)\times\ihom(\olv,\olv)\to\ihom(\olv,\olv),
\]
whose $\Lambda$-component acts on $\lambda\otimes f,\lambda'\otimes g\in\ihom(\olv,\olv)(\Lambda)$ as
\[
\overline{\mu}_\Lambda(\lambda\otimes f,\lambda'\otimes g)=\lambda'\lambda\otimes\mu(f,g)
\]
in accord with our general definition (cf. eq.(\ref{fbar})). Therefore, we can represent the elements of
$\ihom(\olv,\olv)(\Lambda)$ also as matrices, whose entries are now elements of $\Lambda$ of the appropriate
parity. In particular, every element $A\in\ihom(\olv,\olv)(\Lambda)$ can be written as $A=a+c$, where
$a$ is not proportional to any element of $\Lambda$, i.e.,
\[
a=\ihom(\olv,\olv)(\epsilon_\Lambda)(A)
\]
where $\epsilon_\Lambda:\Lambda\to\fieldK$ is the terminal morphism and $c$ is nilpotent. But this shows
that $A$ is invertible if and only if $a$ is invertible, and its inverse is
\begin{equation}
A^{-1}=a^{-1}\sum_{n=0}^\infty(-1)^n(ca^{-1})^n.
\end{equation}
The sum terminates after finitely many terms, since $\Lambda$ is finitely generated. 
Thus the supergroup $\cgl(V)$ is a superdomain, namely
\begin{equation}
\cgl(V)=\ihom(V,V)\big|_{GL(V_{\bar{0}})}.
\end{equation}
The entries of the matrices representing the elements of $\ihom(V,V)=\iend(V)$ are
therefore suited to be used as coordinates of the groups $\cgl(V)$ and their subgroups.

The above results have, of course, been known for a long time and can be obtained without the categorical
framework. But not all groups are that simple. In addition, the above reasoning also implies that
the automorphisms of infinite-dimensional super vector spaces can be extended to Lie supergroups and that
their supermanifold structure is given as a restriction of $\ihom(V,V)$ to $GL(V_{\bar{0}})$.
Another nontrivial example is the diffeomorphism supergroup $\sdiff(\cm)$ of
a supermanifold $\cm$, which will be studied in a subsequent paper \cite{S:Structure}.

\section{Connection to the Berezin-Leites and Rogers-de Witt approaches}

\subsection{Recovery of the ringed space in the finite-dimensional case}

It is instructive to recover the standard ringed space version of a supermanifold from the categorical 
construction. This will, in addition, allow us to relate the categorical and the Berezin-Leites approach
to the Rogers-de Witt approach.
It will also offer a new way to understand superfunctions as actual ``maps'' from somewhere
to somewhere, although these maps will be functor morphisms.
The idea developed here was first described by Molotkov \cite{M:Infinite}. 

One defines an $\olr$-superalgebra $\fr$ in $\catsman$ by setting
\begin{eqnarray}
\fr(\Lambda) &:=& \Lambda\\
\fr(\varphi) &:=& \varphi\quad\textrm{ for }\varphi:\Lambda\to\Lambda'.
\end{eqnarray} 
The $\olr$-superalgebra structure is provided by the $\Lambda_{\bar{0}}$-superalgebra structures on
each $\Lambda$. Note that up to now, we never defined super objects in any of our categories of
superobjects. It was never necessary -- in fact one of the great advantages of the categorical approach is
that one can work with purely even objects. Every
$\fieldK$-supermodule is an ordinary $\olk$-module in $\catsets^\catgr$. In this sense
$\fr$ is ``super super''. The functor $\fr$ is still superrepresentable (and hence indeed an object
of $\catsman$) as we will show now, but the
superalgebra representing it is non-supercommutative.

As an $\olr$-module, we have
\begin{equation}
\fr\cong\olr\oplus\olpi(\olr)\cong\ovlr{1|1}. 
\end{equation}
We want to find a superalgebra structure on $\realR^{1|1}$ which represents the one on $\fr$.
Denoting the standard basis of $\realR^{1|1}$ as $\{1,\theta\}$, we can write
\begin{equation}
\fr(\Lambda)=\Lambda_{\bar{0}}\oplus\Lambda_{\bar{1}}\cong%
(\Lambda_{\bar{0}}\otimes 1)\oplus (\Lambda_{\bar{1}}\otimes\theta).
\end{equation}
Let $\olmu:\fr\times\fr\to\fr$ denote the multiplication in $\fr$ and let $\mu$ denote the
hypothetical multiplication in $\realR^{1|1}$ that we want to determine. Let 
$\lambda_1,\lambda_2\in\Lambda_{\bar{0}}$ be given. We have
\begin{equation}
\olmu_\Lambda(\lambda_1\otimes 1,\lambda_2\otimes 1)=\lambda_2\lambda_1\otimes\mu(1,1).
\end{equation}
Multiplying $\lambda_1,\lambda_2$ within $\fr(\Lambda)$ yields $\lambda_1\lambda_2$, which gets
identified with $\lambda_1\lambda_2\otimes 1$ in $\ovlr{1|1}(\Lambda)$.
Thus we have to require $\mu(1,1)=1$. Likewise, for
$\lambda_1\in\Lambda_{\bar{0}}$ and $\lambda_2\in\Lambda_{\bar{1}}$ we have
\begin{equation}
\olmu_\Lambda(\lambda_1\otimes 1,\lambda_2\otimes\theta)=\lambda_2\lambda_1\otimes\mu(1,\theta).
\end{equation}
This must coincide with $\lambda_1\lambda_2\otimes\theta$, and thus we must have $\mu(1,\theta)=\theta$.
Analogously we find $\mu(\theta,1)=\theta$.

For $\lambda_1,\lambda_2\in\Lambda_{\bar{1}}$, however,
\begin{equation}
\olmu_\Lambda(\lambda_1\otimes\theta,\lambda_2\otimes\theta)=\lambda_2\lambda_1\otimes\mu(\theta,\theta)
\end{equation}
enforces $\mu(\theta,\theta)=-1$.

The super space $\realR^{1|1}$ endowed with this multiplication will be denoted as $\complexC^s$.
As an $\realR$-algebra it is isomorphic to $\complexC$, but as an $\realR$-superalgebra, it is
isomorphic to $\complexC$ with $i$ declared odd. It is non-supercommutative, in particular, every
non-zero element is invertible (which makes it a kind of super analog of a skew field). Note that,
although $\complexC^s$ is non-supercommutative, the $\olr$-superalgebra $\fr$ \emph{is}
supercommutative. Just like passing to their functors of points turns supercommutative algebras
into commutative $\olr$-algebras, the special non-supercommutativity of $\complexC^s$ is weakened
to supercommutativity of its functor of points as an $\olr$-algebra. It would be an interesting
question to study whether this kind of reasoning can be iterated to yield something like 
``super super super'' objects and whether these would possess any geometric interpretation. 
Molotkov in \cite{M:Infinite} proposes a formalism to investigate such
questions, but a conclusive answer has yet to be found.

The reason why we introduced $\fr$ is that we need
a superalgebra in $\catsman$ in order to induce the structure of a superalgebra on certain sets
of morphisms which we want to interpret as the superfunctions on a supermanifold $\cm$. Consider
the set of functor morphisms
\begin{equation}
SC^\infty(\cm):=SC^\infty(\cm,\fr).
\end{equation}
Since $\fr$ is a supercommutative $\olr$-superalgebra, $SC^\infty(\cm)$
is canonically equipped with the structure of a supercommutative $SC^\infty(\cm,\olr)$-su\-per\-al\-ge\-bra.
Moreover, we can embed $\realR\hookrightarrow SC^\infty(\cm)$ as the constant functions $\cm\to\realR$.
More precisely, for any $r\in\realR$, we define a supersmooth morphism $f_r:\cm\to\realR$ by setting
\begin{equation}
(f_r)_\Lambda(m)=r\quad\textrm{ for all }\Lambda\in\catgr,\quad m\in\cm(\Lambda).
\end{equation}
Via this embedding, $SC^\infty(\cm)$ becomes endowed with an $\realR$-superalgebra structure. We
will call $SC^\infty(\cm)$ the algebra of superfunctions on $\cm$.

To see that this construction indeed produces the algebras of superfunctions which are used in the
Berezin-Leites approach, consider the linear supermanifold $\olv$ associated with a finite-dimensional 
super vector space $V$. By Thm.~\ref{thm:skel},
any supersmooth functor morphism $f:\olv\to\fr$ is given by a smooth map $f_0:V_{\bar{0}}\to\realR$ and
by smooth maps
\[
V_{\bar{0}}\to\Sym^n(V_{\bar{1}},\realR^{1|1})\qquad n\geq 1.
\]
The latter maps are just smooth maps
\[
V_{\bar{0}}\to\wedge^n(V_{\bar{1}},\realR)\cong\wedge^n V_{\bar{1}}^*
\]
where $V_{\bar{1}}$ is now considered as an ordinary $\realR$-vector space, $V_{\bar{1}}^*$ is its dual
space, and a map has as parity the exterior degree.
Altogether, the set of these morphisms is therefore isomorphic to the superalgebra
\begin{equation}
C^\infty(V_{\bar{0}})\otimes_\realR\wedge^\bullet(V_{\bar{1}}).
\end{equation}
The superalgebra $SC^\infty(\olv,\fr)$ can then also be written in terms of coordinate maps 
$x_1,\ldots,x_n,\theta_1,\ldots,\theta_m$ for $\olv$, as one usually does.

Now let $\cm$ be a supermanifold and let $M$ be its underlying topological space (i.e., the
topological space underlying the base manifold $\cm(\realR)$). Then we can assign to every open set
$U\subset M$ the $\realR$-superalgebra $SC^\infty(\cm\big|_U)$ and to every inclusion $U'\subset U$
of open sets the obvious restriction map.
This yields a presheaf on $M$ which is easily seen to satisfy the gluing axiom, i.e., it is actually a
sheaf $S(\cm)$. Any morphism $f:\cm\to\cm'$ of supermanifolds
induces, via its associated map $M\to M'$ of the underlying spaces, a morphism of sheaves
$S(f):S(\cm)\to S(\cm')$. Therefore the assignment 
\begin{eqnarray}
\cs:\cm &\mapsto& S(\cm)\\
f &\mapsto& S(f)
\end{eqnarray}
defines a functor from the category
of supermanifolds to the category of topological spaces ringed by supercommutative superalgebras.

Denote by $\catfinsman$ the category of finite-dimensional supermanifolds defined by the categorical
construction and let the category of Berezin-Leites supermanifolds be the supermanifolds obtained
from the standard ringed space construction described in Sect.~\ref{sect:sman}.

\begin{thm}[see also \cite{M:Infinite}]
The functor $\cs$ establishes an equivalence between the category $\catfinsman$ and the category
of Berezin-Leites supermanifolds.
\end{thm}

\begin{proof}
The proof works exactly as in the non-super case. Given a Berezin-Leites supermanifold $\cm=(M,\co_\cm)$
of dimension $m|n$,
its sheaf is defined by the property that every of its stalks $\co_{\cm,x}$ for $x\in M$ is isomorphic
to a superalgebra
\[
\co_{\cm,x}\cong \co_{\fieldK^m,p}\otimes\wedge^\bullet[\theta_1,\ldots,\theta_n]
\]
with odd quantities $\theta_j$, $\co_{\fieldK^m}$ being the sheaf of smooth or holomorphic functions and $p\in\fieldK^m$.
It is therefore isomorphic to a stalk of the sheaf $S(\olv)$ for some super vector space $V$ of dimension
$m|n$. In finite dimensions, i.e., if the stalks are finitely generated, this isomorphism of stalks can
always be extended to a local isomorphism of the corresponding sheaves, because a sheaf is defined
as a local homeomorphism onto its base space. This construction produces then a supermanifold defined in
terms of charts by superdomains in $\olv$, so the functor $\cs$ has a quasi-inverse.
\end{proof}

As is the case in ordinary geometry, this equivalence in general fails to hold if the stalks are not finitely
generated, so in particular in infinite dimensions.

\subsection{The de Witt-Rogers approach}

In the de Witt approach \cite{D:Supermanifolds}, one works over a fixed Grassmann 
algebra $\Lambda$ which is either finitely
generated or the direct limit $\Lambda_\infty$ of the inclusions
\[
\Lambda_{n-1}\hookrightarrow\Lambda_n,\qquad\theta_i\to\theta_i,\quad 1\leq i\leq n-1.
\]

This algebra is then often called the ``supernumbers'', and
it indeed plays the role of numbers in this approach: superfunctions are thought of as functions 
taking values in $\Lambda$, and superdomains of dimension $m|n$ are constructed as open domains in
\[
(\Lambda_{\bar{0}})^m\times(\Lambda_{\bar{1}})^n.
\]
``Open'' refers in this context to the de-Witt topology: a subset $\cu\subset(\Lambda_{\bar{0}})^m\times(\Lambda_{\bar{1}})^n$ 
is called open in this topology if and only if it has the form
\[
\cu=U\times(\Lambda_{\bar{0}}^{nil})^m\times(\Lambda_{\bar{1}})^n,
\]
where $U\subset\fieldK^m$ is an open set in the ordinary sense.

For finitely generated $\Lambda$, super geometry in the de Witt sense coincides with the geometry of the
set of $\Lambda$-points of a functor in $\catsets^\catgr$. 
An open domain of dimension $m|n$ in the de Witt topology is isomorphic to
the set of $\Lambda$-points of an open subfunctor $\cu\subset\olv$, where $\olv$ is a superrepresentable
$\olk$-module. Moreover, the superfunctions on a de Witt supermanifold are just the $\Lambda$-components
$f_\Lambda$ of functor morphisms $f:\cm\to\fr$, as seen in the previous subsection.

To work with a fixed $\Lambda$ is alright \emph{as long as one keeps all constructions functorial
under base change}. This means that when one carries out some construction first for
some Grassmann algebra $\Lambda$ and then for $\Lambda'$, then any morphism $\varphi:\Lambda\to\Lambda'$
must transform the first construction into the second. This excludes, in particular, any attempts to give
the odd dimensions a different topology than the trivial one (cf.~Prop.~\ref{subf}). 
An example for the attempt to avoid artifical side effects
from the choice of $\Lambda$ is the prescription that if one has chosen $\Lambda=\Lambda_n$ (i.e., the
algebra on $n$ generators), the coefficient functions of any superfunctions must be restricted to
$\Lambda_{n-1}$ in order to avoid pathologies in the definition of partial derivatives 
\cite{R:Graded}. Constructions of this sort are completely unnecessary if one sticks to the more
natural requirement of invariance under base change.

The key advantage of the de Witt approach over the ringed space and the categorical approach is
its concreteness. Since the set of $\Lambda$-points of a supermanifold form an actual manifold
and $\fr(\Lambda)=\Lambda$ is just a superalgebra, one can preserve much of the intuitive formalism
of ordinary differential geometry. Superfunctions are still maps from somewhere to somewhere,
and domains are open sets of points in some topology. In most, but not all, practical situations it is sufficient
to work with one set of $\Lambda$-points by choosing a ``generic'' one, which usually just means that $\Lambda$
has to be big enough. Indeed, the vast majority of applications, especially in physics, has been obtained that way.

One just has to keep in mind that both from a mathematical but also from a physical point of view, 
the choice of a fixed
$\Lambda$ can only be of an auxiliary nature. On the side of physics, the choice of a fixed $\Lambda$ and
the use of geometric constructions depending explicitly on it would introduce $\Lambda$ as a fundamental
ingredient of the theories constructed with it. There is no reason to believe that some particular
Grassmann algebra $\Lambda$ plays a special role in nature by providing the odd parameters appearing in
supersymmetric field theories. Moreover, the group $\Aut(\Lambda)$ would assume the illustrous role
of a group of invariances for all physical theories built this way. The only piece of information about
an odd coordinate in a superspace that one is really allowed to fix is that it 
is odd --- anything beyond that
is an arbitrary choice and introduces additional ``fake information''.

\bibliography{complete}
\bibliographystyle{plain}

\noindent
\textbf{email}: sachse@{}mis.mpg.de, \newline
\textbf{Address}: Max Planck Institute for Mathematics in the Sciences, Inselstr. 22, 04103 Leipzig, Germany

\end{document}